\documentclass[hidelinks,10pt,conference,letterpaper]{IEEEtran}

\IEEEoverridecommandlockouts

\usepackage{subcaption}
\usepackage{wrapfig}

\let\labelindent\relax
\usepackage[table,xcdraw]{xcolor}
\usepackage{xspace}
\usepackage{graphicx}
\usepackage{ntheorem}
\usepackage{times,amsmath,epsfig,amssymb,graphicx,amsfonts}
\usepackage{empheq}
\usepackage{graphicx}
\usepackage{psfrag}
\usepackage{fge}
\usepackage{amsmath}
\usepackage{setspace}
\usepackage{color}
\usepackage{amsfonts}   
\usepackage{amssymb}    
\usepackage{mathrsfs}
\usepackage{setspace}
\usepackage{dblfloatfix}
\usepackage{tikz}
\usepackage{tkz-tab}

\usepackage{pgfplotstable}
\usepackage{pgfplots}
\usepackage{algorithm}
\usepackage{algcompatible}
\usepackage{lipsum}
\usepackage{color}
\usepackage{comment}
\usepackage{mathrsfs}
\usepackage{epstopdf}
\usepackage{enumerate}
\usepackage{enumitem}
\usepackage{url}
\usepackage{cite}
\usepackage{tabularx}
\usepackage[pass]{geometry}
\usepackage{chngcntr}

\newcommand{\bdmath}{\begin{dmath}}
\newcommand{\edmath}{\end{dmath}}
\newcommand{\beq}{\begin{equation}}
\newcommand{\eeq}{\end{equation}}
\newcommand{\bdm}{\begin{displaymath}}
\newcommand{\edm}{\end{displaymath}}
\newcommand{\bea}{\begin{eqnarray}}
\newcommand{\eea}{\end{eqnarray}}
\newcommand{\beal}{\beq \begin{array}{lll}}
\newcommand{\eeal}{\end{array} \eeq}
\newcommand{\beas}{\begin{eqnarray*}}
\newcommand{\eeas}{\end{eqnarray*}}
\newcommand{\ba}{\begin{array}}
\newcommand{\ea}{\end{array}}
\newcommand{\bit}{\begin{itemize}}
\newcommand{\eit}{\end{itemize}}
\newcommand{\ben}{\begin{enumerate}}
\newcommand{\een}{\end{enumerate}}

\newcommand{\calA}{{\cal A}}
\newcommand{\calB}{{\cal B}}

\newcommand{\calI}{{\cal I}}

\newcommand{\calM}{{\cal M}}

\newcommand{\calP}{{\cal P}}

\newcommand{\calR}{{\cal R}}
\newcommand{\calS}{{\cal S}}

\newcommand{\calY}{{\cal Y}}
\newcommand{\calV}{{\cal V}}

\newcommand{\calX}{{\cal X}}
\newcommand{\calZ}{{\cal Z}}

\newcommand{\hide}[1]{}

\newcommand{\hiddenText}{{\color{gray} hidden text.}}
\newcommand{\hideWithText}[1]{\hiddenText}

\newcommand{\diag}[1]{\mathrm{diag}\left(#1\right)}

\newcommand{\eye}{{\mathbf I}}

\newcommand{\Real}[1]{ { {\mathbb R}^{#1} } }

\newcommand{\LQG}{LQG\xspace}
\newcommand{\myParagraph}[1]{{\bf #1.}\xspace}

\newcommand{\tlogdet}{{\tt greedy}\xspace}
\newcommand{\tslqg}{{\tt resilient}\xspace}
\newcommand{\toptimal}{{\tt optimal}\xspace}

\newcommand{\validated}[2]{{#2}}

\newcommand{\elem}{{{v}}}

\floatname{algorithm}{Algorithm}

\floatname{algorithm}{Algorithm}

\newtheorem{mydef}{Definition}

\newtheorem{mytheorem}{Theorem}
\newtheorem{mylemma}{Lemma}
\newtheorem{myremark}{Remark}
\newtheorem{mycorollary}{Corollary}

\newtheorem{myproblem}{Problem}

\newcounter{ale}

\newenvironment{liste}{\begin{itemize}}{\end{itemize}}
\newcommand{\aliste}{\begin{liste} \setcounter{ale}{1}}
\newcommand{\zliste}{\end{liste}}

\title{\LARGE \bf Resilient Non-Submodular Maximization over Matroid Constraints}
\author{Vasileios Tzoumas,{$^{1}$} Ali Jadbabaie,{$^{2}$} George J.~Pappas{$^{1}$}
\thanks{$^{1}$The authors are with the Department of Electrical and Systems Engineering, University of Pennsylvania, Philadelphia, PA 19104-6228 USA (email: {\fontsize{8}{8}\selectfont\ttfamily\upshape \{vtzoumas, pappasg\}@seas.upenn.edu}).}
\thanks{$^{2}$The author is with the Institute for Data, Systems and Society, Massachusetts Institute of Technology, Cambridge, MA 02139 USA (email: {\fontsize{8}{8}\selectfont\ttfamily\upshape  jadbabai@mit.edu}).}
\thanks{This work was partially supported  
by the AFOSR Complex Networks program, by the ARL CRA DCIST W911NF-17-2-0181 program, and the Rockefeller Foundation.}
}

\begin{document}
\maketitle

\begin{abstract}
The control and sensing of large-scale systems results in combinatorial problems not only for sensor and actuator placement but also for scheduling or observability/controllability.  Such combinatorial constraints in system design and implementation can be captured using a structure known as matroids.  In particular, the algebraic structure of matroids can be exploited to develop scalable algorithms for sensor and actuator selection, along with quantifiable approximation bounds.  However, in large-scale systems, sensors and actuators may fail or may be \mbox{(cyber-)}attacked.  The objective of this paper is to focus on resilient matroid-constrained problems arising in control and sensing but in the presence of sensor and actuator failures.  In general, resilient matroid-constrained problems are computationally hard. Contrary to the non-resilient case (with no failures), even though they often involve objective functions that are monotone or submodular, no scalable approximation algorithms are known for their solution.  In this paper, we provide the first algorithm, that also has the following properties:  First, it achieves system-wide resiliency, i.e., the algorithm is valid for any number of denial-of-service attacks or failures.  Second, it is scalable, as our algorithm terminates with the same running time as state-of-the-art algorithms for (non-resilient) matroid-constrained optimization.  Third, it provides provable approximation bounds on the system performance, since for monotone objective functions our algorithm guarantees a solution close to the optimal. We quantify our algorithm's approximation performance using a notion of curvature for monotone (not necessarily submodular) set functions. Finally, we support our theoretical analyses with numerical experiments, by considering a control-aware sensor selection scenario, namely, sensing-constrained robot navigation.
\end{abstract}

\section{Introduction}\label{sec:Intro}

Applications in control, robotics, and optimization require system designs in problems such as:
\begin{itemize}
\item (\textit{Control}) \textit{Sensor placement}: In a linear time-invariant system, which few state variables should we measure to ensure observability?~\cite{PEQUITO2017261} 
\item (\textit{Robotics}) \textit{Motion scheduling}: At a team of flying robots, how should we schedule the robots' motions to ensure their capability to track a maximal number of mobile targets?~\cite{tokekar2014multi}
\item (\textit{Optimization}) \textit{Data selection}: Given a flood of driving data collected from the smart-phones of several types of drivers (e.g., truck or commercial vehicle drivers),  which
{few} data should we process from each driver-type to maximize the prediction accuracy of car traffic?~\cite{calinescu2007maximizing}
\end{itemize}
In particular, all above applications motivate the selection of sensors and actuators subject to complex combinatorial constraints, called \textit{matroids}~\cite{fisher1978analysis}, which require the sensors and actuators not only to be a few in number but also to satisfy 
\textit{heterogeneity} restrictions (e.g., few driver data across each driver-type), or \textit{interdependency} restrictions (e.g., system observability).  Other applications in control, robotics, and optimization that involve matroid constraints are:
\begin{itemize}
\item (\textit{Control}) Sparse actuator and sensor scheduling~\cite{summers2016actuator,tzoumas2016minimal,tzoumas2016near,zhang2017kalman,jawaid2015submodularity}; stabilization and voltage control in power grids~\cite{liu2017submodular,liu2018submodular}; and synchronization in complex networks~\cite{clark2017toward};
\item (\textit{Robotics}) Task allocation in collaborative multi-robot systems~\cite{williams2017matroid}; and agile autonomous robot navigation and sparse visual-cue selection~\cite{carlone2016attention}; 
\item (\textit{Optimization}) Sparse signal
recovery and subset column selection~\cite{candes2006stable,boutsidis2009improved,elenberg2016restricted}; and sparse approximation, dictionary and feature selection~\cite{cevher2011greedy,das2011spectral,khanna2017scalable}.
\end{itemize} 

In more detail, all the aforementioned applications~\cite{PEQUITO2017261,tokekar2014multi,calinescu2007maximizing,summers2016actuator,tzoumas2016minimal,tzoumas2016near,zhang2017kalman,jawaid2015submodularity,liu2017submodular,liu2018submodular,clark2017toward,williams2017matroid,carlone2016attention,candes2006stable,boutsidis2009improved,elenberg2016restricted,cevher2011greedy,das2011spectral,khanna2017scalable} require the solution to an optimization problem of the form:
\begin{equation}\label{eq:non_res}
\underset{\!\mathcal{A}\subseteq \calV, \;\calA\in \calI}{\max} \; \;\; f(\mathcal{A}),
\end{equation}
where the set $\calV$ represents the available elements to choose from (e.g., the available sensors); the set $\calI$ represents a matroid constraint (e.g., a cardinality constraint on the number of sensors to be used, and/or a requirement for system observability); and the function~$f$ is a monotone and \textit{possibly} submodular objective function (submodularity is a diminishing returns property).  In particular, $f$ can capture a performance objective (e.g., estimation accuracy).
Notably, the problem in eq.~\eqref{eq:non_res} is combinatorial, and, specifically, is NP-hard~\cite{Feige:1998:TLN:285055.285059}; notwithstanding, approximation algorithms have been proposed {for its solution, such as the greedy~\cite{iyer2013curvature,bian2017guarantees,Feige:1998:TLN:285055.285059,fisher1978analysis,conforti1984curvature}.}

But in all above critical applications, actuators can fail~\cite{willsky1976survey}; sensors can get (cyber-)attacked~\cite{wood2002dos}; and data can get deleted~\cite{mirzasoleiman2017deletion}.  Hence, in such failure-prone and adversarial scenarios,  \textit{resilient} matroid-constrained designs against denial-of-service attacks, failures, or deletions become important. 

In this paper, we formalize for the first time a problem  of \textit{resilient non-submodular maximization}, that goes beyond the traditional problem in eq.~\eqref{eq:non_res}, and guards against attacks, failures, and deletions.  In~particular, we introduce the following resilient re-formulation of the problem in eq.~\eqref{eq:non_res}:
\begin{equation}\label{eq:res}
\underset{\!\mathcal{A}\subseteq \calV, \;\calA\in \calI}{\max} \; \; \underset{\;\mathcal{B}\subseteq \calA, \;\calB\in \calI'}{\min}   \; \;\; f(\mathcal{A}\setminus \mathcal{B}),
\end{equation}
where the set $\calI'$ represents the collection of possible set-removals $\calB$ ---attacks, failures, or deletions--- from $\calA$, each of some specified cardinality. Overall, the problem in eq.~\eqref{eq:res}
maximizes~$f$ despite \textit{worst-case} failures that compromise the maximization in eq.~\eqref{eq:non_res}. Therefore, the problem formulation in eq.~\eqref{eq:res} is suitable in scenarios where there is no prior on the removal mechanism, as well as, in scenarios where protection against worst-case removals is essential, such as in expensive experiment designs, or missions of adversarial-target tracking. 

Particularly, the optimization problem in eq.~\eqref{eq:res} may be interpreted as a $2$-stage perfect information sequential game between two players~\cite[Chapter~4]{myerson2013game}, namely, a ``maximization'' player (designer), and a ``minimization'' player (attacker), where the designer plays first, by selecting $\calA$ to maximize the objective function $f$\!, and, in contrast, the attacker plays second, by selecting $\calB$ to minimize the objective function~$f$\!.  In more detail, \emph{the attacker first observes the selection~$\calA$}, and then, selects $\calB$ such that $\calB$ is a worst-case set removal from~$\calA$. 

In sum, the optimization problem in eq.~\eqref{eq:res} goes beyond traditional (non-resilient) optimization~\cite{Feige:1998:TLN:285055.285059,iyer2013curvature,bian2017guarantees,fisher1978analysis,conforti1984curvature} by proposing \textit{resilient} optimization;
beyond merely cardinality-constrained resilient optimization~\cite{orlin2015robust,tzoumas2017resilient,bogunovic2018robust} by proposing \textit{matroid-constrained} resilient optimization; 
and beyond protection against {non}-adversarial set-removals~\cite{mirzasoleiman2017deletion,kazemi2017deletion} by proposing protection against \textit{worst-case} set-removals.  Hence, the problem in eq.~\eqref{eq:res} aims to protect the complex design of systems, per {heterogeneity} or {interdependency} constraints, against attacks, failures, or deletions, which is a vital objective both for the safety of critical infrastructures, such as power grids~\cite{liu2017submodular,liu2018submodular}, and for the safety of critical missions, such as multi-target surveillance with teams of mobile robots~\cite{tokekar2014multi}.
\medskip

\myParagraph{Contributions} In this paper, we make the contributions:
\begin{itemize}
\item (\textit{Problem}) We formalize the problem of \textit{resilient maximization over matroid constraints} against denial-of-service removals, per eq.~\eqref{eq:res}. This is the first work to formalize, address, and motivate this problem. 
\item (\textit{Solution}) We develop the first algorithm for the problem of resilient maximization over matroid constraints in eq.~\eqref{eq:res}, and prove it enjoys the following properties:
\begin{itemize}
\item \textit{system-wide resiliency}: the algorithm is valid for any number of removals;
\item \textit{minimal running time}: the algorithm terminates with the same running time as state-of-the-art algorithms for (non-resilient) matroid-constrained optimization;
\item \textit{provable approximation performance}: for functions $f$ that are monotone and (possibly) submodular 
---as it holds true in all above applications~\cite{PEQUITO2017261,tokekar2014multi,calinescu2007maximizing,carlone2016attention,summers2016actuator,tzoumas2016minimal,tzoumas2016near,zhang2017kalman,candes2006stable,boutsidis2009improved,elenberg2016restricted,liu2018submodular,jawaid2015submodularity,clark2017toward,cevher2011greedy,das2011spectral,khanna2017scalable,liu2017submodular,williams2017matroid},--- the algorithm ensures a solution close-to-optimal.  

To quantify the algorithm's approximation performance, we use a notion of curvature for monotone (not necessarily submodular) set functions.
\end{itemize}
\item (\textit{Simulations}) We demonstrate the necessity for the resilient re-formulation of the problem in eq.~\eqref{eq:non_res} by conducting numerical experiments in various scenarios of sensing-constrained autonomous robot navigation, varying the number of sensor failures. In addition, via these experiments we demonstrate the benefits of our approach.
\end{itemize}

Overall, the proposed algorithm herein enables the resilient re-formulation and solution of all aforementioned matroid-constrained applications in control, robotics, and optimization~\cite{PEQUITO2017261,tokekar2014multi,calinescu2007maximizing,carlone2016attention,summers2016actuator,tzoumas2016minimal,tzoumas2016near,zhang2017kalman,candes2006stable,boutsidis2009improved,elenberg2016restricted,liu2018submodular,jawaid2015submodularity,clark2017toward,cevher2011greedy,das2011spectral,khanna2017scalable,liu2017submodular,williams2017matroid}; we describe in detail the matroid-constraints involved in all aforementioned application in Section~\ref{sec:problem_statement}.  Moreover, the proposed algorithm enjoys minimal running time, and provable approximation guarantees. 

\medskip

\myParagraph{Organization of the rest of the paper} 
Section~\ref{sec:problem_statement} formulates the problem of resilient maximization over matroid constraints (Problem~\ref{pr:robust_sub_max}), and describes types of matroid constraints in control, robotics, and optimization. Section~\ref{sec:algorithm} presents the first scalable, near-optimal algorithm for Problem~\ref{pr:robust_sub_max}. Section~\ref{sec:performance} presents the main result in this paper, which characterizes the scalability and performance guarantees of the proposed algorithm. Section~\ref{sec:simulations} presents numerical experiments over a control-aware sensor selection scenario. Section~\ref{sec:con} concludes the paper. All proofs are found in the Appendix.

\medskip

\myParagraph{Notation}   
Calligraphic fonts denote sets (e.g., $\calA$).  Given a set $\calA$, then $2^\calA$ denotes the power set of $\calA$; $|\calA|$ denotes $\calA$'s cardinality; given also a set $\calB$, then $\calA\setminus\calB$ denotes the set of elements in $\calA$ that are not in~$\calB$; and the $(\calA,\calB)$ is equivalent to $\calA\cup\calB$. Given a ground set $\mathcal{V}$, a set function $f:2^\mathcal{V}\mapsto \mathbb{R}$, and an element $x\in \mathcal{V}$, the $f(x)$ denotes $f(\{x\})$.

\section{Resilient Non-Submodular Maximization over Matroid Constraints}\label{sec:problem_statement}

We formally define \emph{resilient non-submodular  maximization over matroid constraints}.  
We start with some basic definitions.

\begin{mydef}[Monotonicity]\label{def:mon}
Consider a finite ground set~$\mathcal{V}$. Then, a set function $f:2^\mathcal{V}\mapsto \mathbb{R}$ is \emph{non-decreasing} if and only if for any sets $\mathcal{A}\subseteq \mathcal{A}'\subseteq\mathcal{V}$, it holds $f(\mathcal{A})\leq f(\mathcal{A}')$.
\end{mydef}

\begin{mydef}[Matroid{~\cite[Section~39.1]{schrijver2003combinatorial}}]\label{def:matroid}
Consider a finite ground set~$\mathcal{V}$, and a non-empty collection of subsets of $\calV$, denoted by~$\calI$.  Then, the pair $(\calV,\calI)$ is called a \emph{matroid} if and only if the following conditions hold:
\begin{itemize}
\item for any set $\calX\subseteq \calV$ such that $\calX\in \calI$, and for any set such that $\calZ \subseteq \calX$, it holds $\calZ\in \calI$;
\item for any sets $\calX,\calZ\subseteq \calV$ such that $\calX,\calZ\in \calI$ and $|\calX|<|\calZ|$, it holds that there exists an element $z \in \calZ\setminus \calX$ such that  $\calX\cup\{z\}\in \calI$.
\end{itemize}
\end{mydef}

We next motivate Definition~\ref{def:matroid} by presenting three matroid examples ---uniform, partition, and transversal matroid--- that appear in applications in control, robotics, and optimization.

\myParagraph{Uniform matroid, and applications} 
A matroid $(\calV,\calI)$ is a \emph{uniform matroid} if for a positive integer~$\alpha$ it holds $\calI\equiv\{\calA: \calA\subseteq \calV, |\calA|\leq \alpha\}$.  
Thus, the uniform matroid treats all elements in $\calV$ \textit{uniformly} (that is, as being the same), by only limiting their number in each set that is feasible in~$\calI$.

Applications of the uniform matroid in control, robotics, and optimization, arise when one cannot use an arbitrary number of system elements, e.g., sensors or actuators, to achieve a desired system performance; for example, such sparse element-selection scenarios are necessitated in resource constrained environments of, e.g., limited battery, communication bandwidth, or data processing time~\cite{carlone2016attention}. In more detail, applications of such sparse uniform selection in control, robotics, and optimization include the following: 
\begin{itemize}
\item (\textit{Control}) Actuator and sensor placement, e.g., for system controllability with minimal control effort~\cite{summers2016actuator,tzoumas2016minimal}, and for optimal smoothing or Kalman filtering~\cite{tzoumas2016near,zhang2017kalman};
\item (\textit{Robotics}) Sparse visual-cue selection, e.g., for agile autonomous robot navigation~\cite{carlone2016attention};
\item (\textit{Optimization}) 
Sparse
recovery and column subset selection, e.g., for experiment design~\cite{candes2006stable,boutsidis2009improved,elenberg2016restricted}.
\end{itemize}  

\myParagraph{Partition matroid, and applications} A matroid $(\calV,\calI)$ is a \emph{partition matroid} if for a positive integer $n$, disjoint sets $\calV_1,\ldots,\calV_n$, and positive integers $\alpha_1,\ldots,\alpha_n$, it holds $\calV\equiv \calV_1\cup\cdots\cup\calV_n$ and $\calI\equiv\{\calA: \calA \subseteq \calV,|\calA\cap \calV_i|\leq \alpha_i,\text{ for all } i=1,\ldots,n\}$.  Hence, the partition matroid goes beyond the uniform matroid by allowing for \textit{heterogeneity} in the elements included in each set that is feasible in $\calI$.  We~give two interpretations of the disjoint sets $\calV_1,\ldots,\calV_n$: the first interpretation considers that $\calV_1,\ldots,\calV_n$ correspond to the available elements across $n$ different \textit{types} (buckets) of elements, and correspondingly, the positive integers $\alpha_1,\ldots,\alpha_n$ constrain uniformly the number of elements one can use from each type $1,\ldots,n$ towards a system design goal; the second interpretation considers that $\calV_1,\ldots,\calV_n$ correspond to the available elements across $n$ different \textit{times}, and correspondingly, the positive integers $\alpha_1,\ldots,\alpha_n$ constrain uniformly the number of elements that one can use at each time $1,\ldots,n$.

Applications of the partition matroid in control, robotics, and optimization include all the aforementioned applications in scenarios where heterogeneity in the element-selection enhances the system performance; for example, to guarantee voltage control in power grids, one needs to (possibly) actuate different types of actuators~\cite{liu2018submodular}, and to guarantee active target tracking, one needs to activate different sensors at each time step~\cite{jawaid2015submodularity}. Additional applications of the partition matroid in control and robotics include the following:
\begin{itemize}
\item (\textit{Control}) Synchronization in complex dynamical networks, e.g., for missions of motion coordination~\cite{clark2017toward};
\item (\textit{Robotics}) Robot motion planning, e.g., for multi-target tracking with mobile robots~\cite{tokekar2014multi};
\item (\textit{Optimization}) 
Sparse approximation and feature selection, e.g., for sparse dictionary selection~\cite{cevher2011greedy,das2011spectral,khanna2017scalable}.
\end{itemize}

\myParagraph{Transversal matroid, and applications} A matroid $(\calV,\calI)$ is a \emph{transversal matroid} if for a positive integer~$n$, and a collection of subsets $\calS_1,\ldots,\calS_n$ of $\calV$, it holds~$I$ is the collection of all partial transversals of $(\calS_1,\ldots,\calS_n)$ ---a {partial transversal} is defined as follows: for a finite set~$\calV$, a positive integer~$n$, and a collection of subsets $\calS_1,\ldots,\calS_n$ of $\calV$, a \emph{partial transversal} of $(\calS_1,\ldots,\calS_n)$ is a subset $\calP$ of $\calV$ such that there exist a one-to-one map $\phi: \calP\mapsto \{1,\ldots,n\}$ so that for all $p\in\calP$ it holds $p\in \calS_{\phi(p)}$; i.e., each element in $\calP$ intersects with one ---and only one--- set among the sets $\calS_1,\ldots,\calS_n$.

An application of the transversal matroid in control is that of actuation selection for optimal control performance subject to structural controllability constraints~\cite{PEQUITO2017261}. 

\myParagraph{Additional examples} Other matroid constraints in control, robotics, and optimization are found in the following papers:
\begin{itemize}
\item (\textit{Control}) \cite{liu2017submodular}, for the stabilization of power grids;
\item (\textit{Robotics}) \cite{williams2017matroid}, for task allocation in multi-robot systems;
\item (\textit{Optimization}) \cite{calinescu2007maximizing}, for general task assignments.
\end{itemize}

Given the aforementioned matroid-constrained application-examples, we now define the main problem in this paper.

\begin{myproblem}\label{pr:robust_sub_max} 
\emph{\textbf{(Resilient Non-Submodular Maximization over Matroid Constraints)}}
Consider the problem parameters: 
\begin{itemize}
\item a matroid $(\calV,\calI)$;
\item a matroid $(\calV,\calI')$ such that: $(\calV,\calI')$ is either a uniform matroid, or a partition matroid $(\calV,\calI')$ with the same partition as $(\calV,\calI)$;
\item a non-decreasing set function $f:2^{\mathcal{V}} \mapsto \mathbb{R}$
such that (without loss of generality) it holds $f(\emptyset)=0$, and for any set $\calA\subseteq \calV$, it also holds $f(\calA)\geq 0$.
\end{itemize}  

The problem of \emph{resilient non-submodular maximization over matroid constraints} is to 
maximize the function~$f$ by selecting a set $\calA\subseteq \calV$ such that $\calA\in \calI$, and accounting for any worst-case set removal $\calB\subseteq \calA$ from $\calA$ such that $\calB\in \calI'$\!\!.
Formally:\footnote{Given a matroid $(\calV,\calI')$, and any subset $\calA\subseteq \calV$, then, the $(\calA,\{\calB: \mathcal{B}\subseteq \calA, \calB\in \calI'\})$ is also a matroid~\cite[Section~39.3]{schrijver2003combinatorial}.}
\vspace*{1mm}
\begin{equation*}
\underset{\mathcal{A}\subseteq \calV,\; \calA\in \calI}{\max} \; \; \underset{\;\mathcal{B}\subseteq \calA,\; \calB\in \calI'}{\min}   \; \;\; f(\mathcal{A}\setminus \mathcal{B}).
\end{equation*}
\end{myproblem}
\medskip

As we mentioned in this paper's Introduction, Problem~\ref{alg:rob_sub_max} may be interpreted as a $2$-stage perfect information sequential game between two players~\cite[Chapter~4]{myerson2013game}, namely, a ``maximization'' player, and a ``minimization'' player, 
where the ``maximization'' player plays first by selecting the set~$\calA$, and, then, \emph{the ``minimization'' player observes $\calA$}, and plays second by selecting a worst-case set removal $\calB$ from $\calA$.

In sum, Problem~\ref{pr:robust_sub_max} aims to guard all the aforementioned applications in control, robotics, and optimization~\cite{PEQUITO2017261,tokekar2014multi,calinescu2007maximizing,carlone2016attention, summers2016actuator,tzoumas2016minimal,tzoumas2016near,zhang2017kalman,candes2006stable,boutsidis2009improved,elenberg2016restricted,liu2018submodular,jawaid2015submodularity,clark2017toward,cevher2011greedy,das2011spectral,khanna2017scalable,liu2017submodular,williams2017matroid} against sensor and actuator attacks or failures, by proposing their resilient re-formulation, since all involve the maximization of non-decreasing functions subject to matroid constrains. 
For example, we discuss the resilient re-formulation of two among the aforementioned applications~\cite{PEQUITO2017261,tokekar2014multi,calinescu2007maximizing,carlone2016attention, summers2016actuator,tzoumas2016minimal,tzoumas2016near,zhang2017kalman,candes2006stable,boutsidis2009improved,elenberg2016restricted,liu2018submodular,jawaid2015submodularity,clark2017toward,cevher2011greedy,das2011spectral,khanna2017scalable,liu2017submodular,williams2017matroid}:
\setcounter{paragraph}{0}
\paragraph{Actuator placement for minimal control effort~\cite{summers2016actuator,tzoumas2016minimal}} Given a dynamical system, the design objective is to select a few actuators to place in the system to achieve controllability with minimal control effort~\cite{tzoumas2016minimal}.  In particular, the actuator-selection framework is as follows: given a set $\calV$ of available actuators to choose from, then, up to $\alpha$ actuators can be placed in the system.  In more detail, the aforementioned actuator-selection problem can be captured by a uniform matroid $(\calV,\calI)$ where $\calI\triangleq\{\calA:\calA\in \calV,|\calA|\leq \alpha\}$. However, in the case of a failure-prone environment where up to~$\beta$ actuators may fail, then a resilient re-formulation of the aforementioned problem formulation is necessary: Problem~\ref{pr:robust_sub_max} suggests that such a resilient re-formulation can be achieved by modelling any set of $\beta$ actuator-failures in $\calA$ by a set $\calB$ in the uniform matroid on $\calA$ where $\calB\subseteq \calA$ and $|\calB|\leq \beta$.
\paragraph{Multi-target coverage with mobile robots~\cite{tokekar2014multi}} A number of adversarial targets are deployed in the environment, and a team of mobile robots~$\calR$ is tasked to cover them using on-board cameras. To this end, at each time step the robots in~$\calR$ need to jointly choose their motion. In particular, the movement-selection framework is as follows: given a finite set of possible moves~$\calM_i$ for each robot $i\in\calR$, then, at each time step each robot selects a move to make so that the team $\calR$ covers collectively as many targets as possible.  In more detail, since each robot in~$\calR$ can make only one move per time, if we denote by~$\calA$ the set of moves to be made by each robot in~$\calR$, then the aforementioned movement-selection problem can be captured by a partition matroid $(\calV,\calI)$ such that $\calV=\cup_{i\in\calR} \calM_i$ and $\calI=\{\calA: \calA\subseteq \calV, |\calM_i\cap\calA|\leq 1, \text{ for all } i\in\calR\}$~\cite{tokekar2014multi}. However, in the case of an adversarial scenario where the targets can attack up to~$\beta$ robots, then a resilient re-formulation of the aforementioned problem formulation is necessary: Problem~\ref{pr:robust_sub_max} suggests that such a resilient re-formulation can be achieved by modelling any set of $\beta$ attacks to the robots in $\calR$ by a set $\calB$ in the uniform matroid on $\calS$ where $\calB\subseteq \calS$ and $|\calB|\leq \beta$. 


\section{Algorithm for Problem~\ref{pr:robust_sub_max}} \label{sec:algorithm}

We present the first scalable algorithm for Problem~\ref{pr:robust_sub_max}. 
The pseudo-code of the algorithm is described in Algorithm~\ref{alg:rob_sub_max}. 

\subsection{Intuition behind Algorithm~\ref{alg:rob_sub_max}}\label{subsec:intuition}

The goal of Problem~\ref{pr:robust_sub_max} is to ensure a maximal value for an objective function $f$ through a single maximization step, despite compromises to the solution of the maximization step.  In~particular, Problem~\ref{pr:robust_sub_max} aims to select a set $\calA$ towards a maximal value of $f$\!, 
despite that $\calA$ is later compromised by a worst-case set removal $\calB$, resulting to $f$ being finally evaluated at the set $\calA\setminus \calB$ instead of the set $\calA$.
In~this~context, Algorithm~\ref{alg:rob_sub_max} aims to fulfil the goal of Problem~\ref{pr:robust_sub_max} by constructing the set $\calA$ as the union of two sets, namely, the $\calA_{1}$ and $\calA_{2}$ (line~\ref{line:selection} of Algorithm~\ref{alg:rob_sub_max}), whose role we describe in more detail below:
\setcounter{paragraph}{0} 
\paragraph{Set $\calA_{1}$ approximates worst-case set removal from~$\calA$}  Algorithm~\ref{alg:rob_sub_max} aims with the set $\calA_{1}$  to capture a worst-case set-removal of elements ---per the matroid $(\calV,\calI')$--- from the elements Algorithm~\ref{alg:rob_sub_max} is going to select in the set~$\calA$; equivalently, the set~$\calA_{1}$ is aimed to act as a ``bait'' to an attacker that selects to remove the \textit{best} set of elements from~$\calA$ per the matroid $(\calV,\calI')$ (\textit{best} with respect to the elements' contribution towards the goal of Problem~\ref{pr:robust_sub_max}). However, the problem of selecting the \textit{best} elements in~$\calV$ per a matroid constraint is a combinatorial and, in general, intractable problem~\cite{Feige:1998:TLN:285055.285059}. 
For this reason, Algorithm~\ref{alg:rob_sub_max} aims to \textit{approximate} the best set of elements in $\calI'$\!\!, by letting $\calA_{1}$ be the set of elements with the largest marginal contributions to the value of the objective function~$f$ (lines~\ref{line:begin_while_1}-\ref{line:end_while_1} of Algorithm~\ref{alg:rob_sub_max}). In addition, since  per Problem~\ref{pr:robust_sub_max} the set $\calA$ needs to be in the matroid $(\calV,\calI)$, Algorithm~\ref{alg:rob_sub_max} constructs $\calA_1$ so that not only it is $\calA_1\in\calI'$\!\!, as we described before, but so that it also is $\calA_1\in \calI$ (lines~\ref{line:begin_if_1}-\ref{line:end_if_1} of Algorithm~\ref{alg:rob_sub_max}).

\paragraph{Set $\calA_{2}$ is such that the set $\calA_{1}\cup \calA_{2}$ approximates optimal solution to Problem~\ref{pr:robust_sub_max}}
Assuming that~$\calA_{1}$ is the set that is going to be removed from Algorithm~\ref{alg:rob_sub_max}'s set selection~$\calA$,
Algorithm~\ref{alg:rob_sub_max} needs to select a set of elements $\calA_{2}$ to complete the construction of~$\calA$ so that~$\calA=\calA_1\cup\calA_2$ is in the matroid $(\calV,\calI)$, per Problem~\ref{pr:robust_sub_max}. In~particular, for $\calA=\calA_{1}\cup \calA_{2}$ to be an optimal solution to Problem~\ref{pr:robust_sub_max} (assuming the removal of~$\calA_{1}$ from $\calA$), Algorithm~\ref{alg:rob_sub_max} needs to select $\calA_{2}$ as a \textit{best} set of elements from $\calV\setminus\calA_{1}$ subject to the constraint that $\calA_1\cup\calA_2$ is in $(\calV,\calI)$ (lines~\ref{line:begin_if_2}-\ref{line:end_if_2} of Algorithm~\ref{alg:rob_sub_max}).  
Nevertheless, the problem of selecting a \textit{best} set 
of elements subject to such a constraint is a combinatorial and, in~general, intractable problem~\cite{Feige:1998:TLN:285055.285059}.  Hence, Algorithm~\ref{alg:rob_sub_max} aims to \textit{approximate} such a best set, 
using the greedy procedure in the lines~\ref{line:begin_while_2}-\ref{line:end_while_2} of   Algorithm~\ref{alg:rob_sub_max}.  

Overall, Algorithm~\ref{alg:rob_sub_max} constructs the sets $\calA_{1}$ and $\calA_{2}$ to approximate with their union $\calA$ an optimal solution to Problem~\ref{pr:robust_sub_max}.

We next describe the steps in Algorithm~\ref{alg:rob_sub_max} in more detail.

\begin{algorithm}[t]
\caption{Scalable algorithm for Problem~\ref{pr:robust_sub_max}.}
\begin{algorithmic}[1]
\REQUIRE Per Problem~\ref{pr:robust_sub_max}, Algorithm~\ref{alg:rob_sub_max} receives the inputs:
\begin{itemize}
\item a matroid $(\calV,\calI)$;
\item an either uniform or partition matroid $(\calV,\calI')$;
\item a non-decreasing set function $f:2^{\mathcal{V}} \mapsto \mathbb{R}$
such that it is $f(\emptyset)=0$, and for any set $\calA\subseteq \calV$, it also is $f(\calA)\geq 0$.
\end{itemize}
\ENSURE Set $\mathcal{A}$.
\medskip

\STATE $\mathcal{A}_{1}\leftarrow\emptyset$;~~~$\mathcal{R}_{1}\leftarrow\emptyset$;~~~$\mathcal{A}_{2}\leftarrow\emptyset$;~~~$\mathcal{R}_{2}\leftarrow\emptyset$;\label{line:initiliaze}
\WHILE {$\mathcal{R}_{1}\neq \calV$}\label{line:begin_while_1}
\STATE $x\in \arg\max_{y \in \calV\setminus\mathcal{R}_{1}} f(y)$;\label{line:select_element_bait}
\IF {$\mathcal{A}_{1}\cup\{x\}\in \calI$ and $\mathcal{A}_{1}\cup\{x\}\in \calI'$\!} \label{line:begin_if_1}
\STATE $\mathcal{A}_{1}\leftarrow\mathcal{A}_{1}\cup\{x\}$;\label{line:build_of_bait}
\ENDIF \label{line:end_if_1}
\STATE {$\mathcal{R}_{1}\leftarrow \mathcal{R}_{1}\cup \{x\}$}; \label{line:increase_removed_set_1}
\ENDWHILE \label{line:end_while_1}
\WHILE {$\mathcal{R}_{2}\neq \calV\setminus \calA_1$} \label{line:begin_while_2} 
\STATE $x\in \arg\max_{y \in \mathcal{V}\setminus (\mathcal{A}_{1}\cup\mathcal{R}_{2})}f(\calA_2\cup \{y\})$; \label{line:greedy_selection}
\IF {$\mathcal{A}_{1}\cup\mathcal{A}_{2}\cup\{x\}\in \calI$} \label{line:begin_if_2}
\STATE $\mathcal{A}_{2}\leftarrow\mathcal{A}_{2}\cup\{x\}$;\label{line:build_of_greedy}
\ENDIF \label{line:end_if_2}
\STATE {$\mathcal{R}_{2}\leftarrow \mathcal{R}_{2}\cup \{x\}$}; \label{line:increase_removed_set_2}
\ENDWHILE \label{line:end_while_2}
\STATE $\mathcal{A}\leftarrow \mathcal{A}_{1} \cup \mathcal{A}_{2}$; \label{line:selection}
\end{algorithmic}\label{alg:rob_sub_max}
\end{algorithm}

\subsection{Description of steps in Algorithm~\ref{alg:rob_sub_max}}

Algorithm~\ref{alg:rob_sub_max} executes four steps:
\setcounter{paragraph}{0}
\paragraph{Initialization (line~\ref{line:initiliaze} of Algorithm~\ref{alg:rob_sub_max})} Algorithm~\ref{alg:rob_sub_max} defines four auxiliary sets, namely, the $\calA_{1}$, $\calR_1$, $\calA_2$, and $\calR_2$, and initializes each of them with the empty set (line~\ref{line:initiliaze} of Algorithm~\ref{alg:rob_sub_max}). \textit{The~purpose of $\calA_{1}$ and $\calA_2$} is to construct the set $\calA$, which is the set Algorithm~\ref{alg:rob_sub_max} selects as a solution to Problem~\ref{pr:robust_sub_max}; in particular, 
the union of $\calA_{1}$ and of~$\calA_2$ constructs~$\calA$ by the end of Algorithm~\ref{alg:rob_sub_max} (line~\ref{line:selection} of Algorithm~\ref{alg:rob_sub_max}).  \textit{The~purpose of $\calR_{1}$ and of $\calR_2$} is to support the construction of $\calA_{1}$ and $\calA_2$, respectively; in particular, during the construction of $\calA_1$, Algorithm~\ref{alg:rob_sub_max} stores in $\calR_1$ the elements of $\calV$ that have either been included already or cannot be included in $\calA_1$ (line~\ref{line:increase_removed_set_1} of Algorithm~\ref{alg:rob_sub_max}), and that way, Algorithm~\ref{alg:rob_sub_max} keeps track of which elements remain to be checked whether they could be added in $\calA_1$ (line~\ref{line:build_of_bait} of Algorithm~\ref{alg:rob_sub_max}). Similarly,  during the construction of $\calA_2$, Algorithm~\ref{alg:rob_sub_max} stores in $\calR_2$ the elements of $\calV\setminus\calA_1$ that have either been included already or cannot be included in $\calA_2$ (line~\ref{line:increase_removed_set_2} of Algorithm~\ref{alg:rob_sub_max}), and that way, Algorithm~\ref{alg:rob_sub_max} keeps track of which elements remain to be checked whether they could be added in $\calA_2$ (line~\ref{line:build_of_greedy} of Algorithm~\ref{alg:rob_sub_max}).

\paragraph{Construction of set $\calA_{1}$ (lines~\ref{line:begin_while_1}-\ref{line:end_while_1} of Algorithm~\ref{alg:rob_sub_max})} Algorithm~\ref{alg:rob_sub_max} constructs the set $\calA_{1}$ sequentially ---by adding one element  at a time from $\calV$ to $\calA_{1}$, over a sequence of multiple time-steps--- such that $\calA_{1}$ is contained in both the matroid $(\calV,\calI)$ and the matroid $(\calV,\calI')$ (line~\ref{line:begin_if_1} of Algorithm~\ref{alg:rob_sub_max}), and such that each element $v\in \calV$ that is chosen to be added in~$\calA_1$ achieves the highest marginal value of $f(v)$ among all the elements in $\calV$ that have not been yet added in~$\calA_1$ and can be added in $\calA_1$ (line~\ref{line:build_of_bait} of Algorithm~\ref{alg:rob_sub_max}).

\paragraph{Construction of set $\calA_{2}$ (lines~\ref{line:begin_while_2}-\ref{line:end_while_2} of Algorithm~\ref{alg:rob_sub_max})} Algorithm~\ref{alg:rob_sub_max} constructs the set $\calA_{2}$ sequentially, by picking greedily elements from the set $\calV_t\setminus \calA_{1}$ such that $\calA_1\cup\calA_2$ is contained in the matroid $(\calV,\calI)$.
Specifically, the greedy procedure in Algorithm~\ref{alg:rob_sub_max}'s ``while loop''  (lines~\ref{line:begin_while_2}-\ref{line:end_while_2} of Algorithm~\ref{alg:rob_sub_max}) selects an element $y\in\mathcal{V}\setminus (\mathcal{A}_{1}\cup\calR_2)$ to add in $\calA_{2}$ only if $y$ maximizes the value of  $f(\calA_2\cup \{y\})$, where the set~$\calR_2$ stores the elements that either have already been added to $\calA_2$ or have been considered to be added to $\calA_2$ but they were not because the resultant set $\calA_1\cup\calA_2$ would not be in the matroid $(\calV,\calI)$.

\paragraph{Construction of set $\calA$ (line~\ref{line:selection} of Algorithm~\ref{alg:rob_sub_max})}
Algorithm~\ref{alg:rob_sub_max} constructs the set $\calA$ as the union of the previously constructed sets $\calA_{1}$ and~$\calA_2$ \mbox{(lines~\ref{line:selection} of Algorithm~\ref{alg:rob_sub_max}).}  

In sum, Algorithm~\ref{alg:rob_sub_max} proposes a set $\calA$ as solution to Problem~\ref{pr:robust_sub_max}, and in particular, Algorithm~\ref{alg:rob_sub_max} constructs the set $\calA$ so it can withstand any compromising set removal from~it.

\section{Performance Guarantees for Algorithm~\ref{alg:rob_sub_max}}\label{sec:performance}

We quantify Algorithm~\ref{alg:rob_sub_max}'s performance, by bounding its running time, and its approximation performance. To this end, we use the following two notions of curvature for set functions, as well as, a notion of rank for a matroid.

\subsection{Curvature and total curvature of non-decreasing functions}\label{sec:total_curvature}
 
We present the notions of \emph{curvature} and of \emph{total curvature} for non-decreasing set functions.  We start by describing the notions of \textit{modularity} and 
 \textit{submodularity} for set functions.

\begin{mydef}[Modularity]\label{def:modular}
Consider any finite set~$\mathcal{V}$.  The set function $g:2^\mathcal{V}\mapsto \mathbb{R}$ is modular if and only if for any set $\mathcal{A}\subseteq \mathcal{V}$, it holds $g(\mathcal{A})=\sum_{\elem\in \mathcal{A}}g(\elem)$.
\end{mydef}

In words, a set function $g:2^\mathcal{V}\mapsto \mathbb{R}$ is modular if through~$g$ all elements in $\mathcal{V}$ cannot substitute each other; in particular, Definition~\ref{def:modular} of modularity implies that for any set $\mathcal{A}\subseteq\mathcal{V}$, and for any element $\elem\in \mathcal{V}\setminus\mathcal{A}$, it holds $g(\{\elem\}\cup\mathcal{A})-g(\mathcal{A})= g(\elem)$.

\begin{mydef}[Submodularity~{\cite[Proposition 2.1]{nemhauser78analysis}}]\label{def:sub}
Consider any finite set $\calV$.  Then, the set function $g:2^\calV\mapsto \mathbb{R}$ is \emph{submodular} if and only if
for any sets $\mathcal{A}\subseteq \validated{\mathcal{A}'}{\mathcal{B}}\subseteq\calV$, and any element $\elem\in \calV$, it \validated{is}{holds}  
$g(\mathcal{A}\cup \{\elem\})\!-\!g(\mathcal{A})\geq g(\validated{\mathcal{A}'}{\mathcal{B}}\cup \{\elem\})\!-\!g(\validated{\mathcal{A}'}{\mathcal{B}})$.
\end{mydef}

Definition~\ref{def:sub} implies that a set function $g:2^\calV\mapsto \mathbb{R}$ is submodular if and only if it satisfies a diminishing returns property where
for any set $\mathcal{A}\subseteq \mathcal{V}$, and for any element $\elem\in \mathcal{V}$, the marginal gain $g(\mathcal{A}\cup \{\elem\})-g(\mathcal{A})$ is~non-increasing. 
In contrast to modularity, submodularity implies that the elements in $\mathcal{V}$ \emph{can} substitute each other, since Definition~\ref{def:sub} of submodularity implies the inequality $g(\{\elem\}\cup\mathcal{A})-g(\mathcal{A})\leq g(\elem)$; that is, in the presence of the set $\mathcal{A}$, the element $\elem$ may lose part of its contribution to the  value of $g(\{x\}\cup\mathcal{A})$.

\begin{mydef}\label{def:curvature}
\emph{\textbf{(Curvature of monotone submodular functions~\cite{conforti1984curvature})}}
Consider a finite set $\mathcal{V}$, and a non-decreasing submodular set function $g:2^\mathcal{V}\mapsto\mathbb{R}$ such that (without loss of generality) for any element $\elem \in \mathcal{V}$, it is  $g(\elem)\neq 0$.  Then, the curvature of $g$ is defined as follows: \begin{equation}\label{eq:curvature}
\kappa_g\triangleq 1-\min_{\elem\in\mathcal{V}}\frac{g(\mathcal{V})-g(\mathcal{V}\setminus\{\elem\})}{g(\elem)}.
\end{equation}
\end{mydef}

Definition~\ref{def:curvature} of curvature implies that for any non-decreasing submodular set function $g:2^\mathcal{V}\mapsto\mathbb{R}$, it holds $0 \leq \kappa_g \leq 1$.  In particular, the value of $\kappa_g$ measures how far~$g$ is from modularity, as we explain next: if $\kappa_g=0$, then for all elements $v\in\mathcal{V}$, it holds $g(\mathcal{V})-g(\mathcal{V}\setminus\{v\})=g(v)$, that is, $g$ is modular. In~contrast, if $\kappa_g=1$, then there exist an element $v\in\mathcal{V}$ such that $g(\mathcal{V})=g(\mathcal{V}\setminus\{v\})$, that is, in the presence of $\mathcal{V}\setminus\{v\}$, $v$~loses all its contribution to the value of $g(\mathcal{V})$.

\begin{mydef}\label{def:total_curvature}
\emph{\textbf{(Total curvature of non-decreasing functions~\cite[Section~8]{sviridenko2017optimal})}}
Consider a finite set $\mathcal{V}$, and a monotone set function $g:2^\mathcal{V}\mapsto\mathbb{R}$.  Them, the total curvature of $g$ is defined as follows: 
\begin{equation}\label{eq:total_curvature}
c_g\triangleq 1-\min_{v\in\mathcal{V}}\min_{\mathcal{A}, \mathcal{B}\subseteq \mathcal{V}\setminus \{v\}}\frac{g(\{v\}\cup\mathcal{A})-g(\mathcal{A})}{g(\{v\}\cup\mathcal{B})-g(\mathcal{B})}.
\end{equation}
\end{mydef}

Definition~\ref{def:total_curvature} of total curvature implies that for any non-decreasing set function  $g:2^\mathcal{V}\mapsto\mathbb{R}$, it holds $0 \leq c_g \leq 1$. To connect the notion of total curvature with that of curvature, we note that when the function $g$ is non-decreasing and submodular, then the two notions coincide, i.e., it holds $c_g=\kappa_g$; the reason is that if $g$ is non-decreasing and submodular, then the inner minimum in eq.~\eqref{eq:total_curvature} is attained for $\calA=\calB\setminus\{v\}$ and $\calB=\emptyset$.  
In addition, to connect the above notion of total curvature with the notion of modularity, we note that 
if $c_g=0$, then $g$ is modular, since eq.~\eqref{eq:total_curvature} implies that for any elements $\elem\in\calV$, and for any sets $\calA,\calB\subseteq \calV\setminus \{\elem\}$, it holds: 
\begin{equation}\label{eq:ineq_total_curvature}
(1-c_g)\left[g(\{\elem\}\cup\calB)-g(\calB)\right]\leq g(\{\elem\}\cup\calA)-g(\calA),
\end{equation}
which for $c_g=0$ implies the modularity of $g$. Finally, 
to connect the above notion of total curvature with the notion of monotonicity, we mention that 
if $c_g=1$, then eq.~\eqref{eq:ineq_total_curvature} implies that $g$ is merely non-decreasing  (as it is already assumed by the Definition~\ref{def:total_curvature} of total curvature).

\subsection{Rank of a matroid}\label{subsec:rank_matroid}

We present a notion of rank for a matroid.
 
\begin{mydef}[Rank of a matroid{~\cite[Section~39.1]{schrijver2003combinatorial}}]\label{def:rank_matroid}
Consider a matroid $(\calV,\calI)$. Then, the rank of $(\calV,\calI)$ is the number equal to the cardinality of the set $\calX\in\calI$ with the maximum cardinality among all sets in $\calI$.
\end{mydef}

For example, per the discussions in Section~\ref{sec:problem_statement}, for a uniform matroid $(\calV,\calI)$ of the form $\calI\equiv\{\calA: \calA\subseteq \calV, |\calA|\leq \alpha\}$, the rank is equal to $\alpha$; and for a partition matroid $(\calV,\calI)$ of the form  $\calV\equiv \calV_1\cup\cdots\cup\calV_n$ and $\calI\equiv\{\calA: \calA \subseteq \calV,|\calA\cap \calV_i|\leq \alpha_i,\text{ for all } i=1,\ldots,n\}$, the rank is equal to $\alpha_1+\ldots+\alpha_n$.

\subsection{Performance analysis for Algorithm~\ref{alg:rob_sub_max}}

We quantify Algorithm~\ref{alg:rob_sub_max}'s approximation performance, as well as, its running time per maximization step in Problem~\ref{pr:robust_sub_max}.

\begin{mytheorem}[Performance of Algorithm~\ref{alg:rob_sub_max}]\label{th:alg_rob_sub_max_performance}
Consider an instance of Problem~\ref{pr:robust_sub_max}, the notation therein, the notation in Algorithm~\ref{alg:rob_sub_max}, and the definitions:
\begin{itemize}
\item let the number $f^\star$ be the (optimal) value to Problem~\ref{pr:robust_sub_max};
\item given a set $\mathcal{A}$ as solution to  Problem~\ref{pr:robust_sub_max}, let  $\mathcal{B}^\star(\mathcal{A})$ be an optimal (worst-case) set removal from $\mathcal{A}$, per Problem~\ref{pr:robust_sub_max}, that is: 
$\mathcal{B}^\star(\mathcal{A})\in\arg\underset{\;\mathcal{B}\subseteq \calA, \calB\in \calI'(\calA)}{\min}   \; \;\; f(\mathcal{A}\setminus \mathcal{B})$;

\item let the numbers $\alpha$ and $\beta$ be such that $\alpha$ is the rank of the matroid $(\calV,\calI)$; and $\beta$ is the rank of the matroid $(\calV,\calI')$; 

\item define $h(\alpha,\beta)\triangleq \max [1/(1+\alpha), 1/(\alpha-\beta)]$.\footnote{A plot of $h(\alpha,\beta)$ is found in Fig.~\ref{fig:bounds2}.}
\end{itemize}

The performance of Algorithm~\ref{alg:rob_sub_max} is bounded as follows:

\begin{enumerate}[leftmargin=*]
\item \emph{(Approximation performance)}~Algorithm~\ref{alg:rob_sub_max} returns a set $\calA$ such that $\mathcal{A}\subseteq \calV$, $\mathcal{A}\in \calI$, and:
\begin{itemize}
\item if $f$ is \emph{non-decreasing} and \emph{submodular}, and:
\begin{itemize}
\item if $(\calV,\calI)$ is a \emph{uniform} matroid, then:
\begin{equation}\label{ineq:bound_sub_uniform}
\frac{f(\mathcal{A}\setminus \mathcal{B}^\star(\calA))}{f^\star}\geq
\frac{\max\left[1-\kappa_f,h(\alpha,\beta)\right]}{\kappa_f}(1-e^{-\kappa_f});
 \end{equation}
\item if $(\calV,\calI)$ is \emph{any} matroid, then:
\begin{equation}\label{ineq:bound_sub}
\frac{f(\mathcal{A}\setminus \mathcal{B}^\star(\calA))}{f^\star}\geq 
\frac{\max\left[1-\kappa_f,h(\alpha,\beta)\right]}{1+\kappa_f};
\end{equation}
\end{itemize} 

\item if $f$ is \emph{non-decreasing}, then:
\begin{equation}\label{ineq:bound_non_sub}
\frac{f(\mathcal{A}\setminus \mathcal{B}^\star(\calA))}{f^\star}\geq  (1-c_f)^3\!\!.
\end{equation}
\end{itemize}

\item \emph{(Running time)}~Algorithm~\ref{alg:rob_sub_max} constructs the set $\calA$ as a solutions to Problem~\ref{pr:robust_sub_max} with $O(|\mathcal{V}|^2)$ evaluations of $f$\!. 
\end{enumerate}
\end{mytheorem}

\myParagraph{Provable approximation performance}
Theorem~\ref{th:alg_rob_sub_max_performance}  implies on the approximation performance of Algorithm~\ref{alg:rob_sub_max}:
\setcounter{paragraph}{0}
\paragraph{{Near-optimality}} Both for any monotone submodular objective functions $f$\!, and for any merely monotone objective functions~$f$ with total curvature $c_f<1$, Algorithm~\ref{alg:rob_sub_max} guarantees a value for Problem~\ref{pr:robust_sub_max} finitely close to the optimal.  In~particular,  per ineq.~\eqref{ineq:bound_sub_uniform} and ineq.~\eqref{ineq:bound_sub} (\textit{case of submodular functions}), the~approximation factor of Algorithm~\ref{alg:rob_sub_max} 
is bounded by $\frac{h_f(\alpha,\beta)}{\kappa_f}(1-e^{-\kappa_f})$ and $\frac{h_f(\alpha,\beta)}{1+\kappa_f}$, respectively, which for any finite  number~$\alpha$ are both
non-zero (see also Fig.~\ref{fig:bounds2}); in addition, per ineq.~\eqref{ineq:bound_sub_uniform} and ineq.~\eqref{ineq:bound_sub}, the~approximation factor of Algorithm~\ref{alg:rob_sub_max} is also bounded by $\frac{1-\kappa_f}{\kappa_f}(1-e^{-\kappa_f})$  and $\frac{1-\kappa_f}{1+\kappa_f}$, respectively, which are also non-zero for any monotone submodular function~$f$ with $\kappa_f<1$ (see also Fig.~\ref{fig:bounds}).
Similarly, per ineq.~\eqref{ineq:bound_non_sub} (\textit{case of monotone functions}), the approximation factor of Algorithm~\ref{alg:rob_sub_max} is bounded by $(1-c_f)^3$\!, which is non-zero for any monotone function~$f$ with $c_f<1$ ---notably, although it is known for the problem of (non-resilient) set function maximization that the approximation bound $(1-c_f)$ is tight~\cite[Theorem~8.6]{sviridenko2017optimal}, the tightness of the bound $(1-c_f)^3$ in ineq.~\eqref{ineq:bound_non_sub} for Problem~\ref{pr:robust_sub_max} is an open problem.

We discuss classes of functions $f$ with curvatures $\kappa_f<1$ or $c_f<1$, along with relevant applications, in the remark below.

\begin{myremark}\emph{\textbf{(Classes of functions $f$ with $\kappa_f<1$ or $c_f<1$, and applications)}}
\emph{Classes of functions $f$ with $\kappa_f<1$} are the concave over modular functions~\cite[Section~2.1]{iyer2013curvature}, and the $\log\det$ of positive-definite matrices~\cite{sviridenko2015optimal,sharma2015greedy}. \emph{Classes of functions $f$ with $c_f<1$} are support selection functions~\cite{elenberg2016restricted}, 
and estimation error metrics such as the average minimum square error of the Kalman filter~\cite[Theorem~4]{tzoumas2018codesign}.

 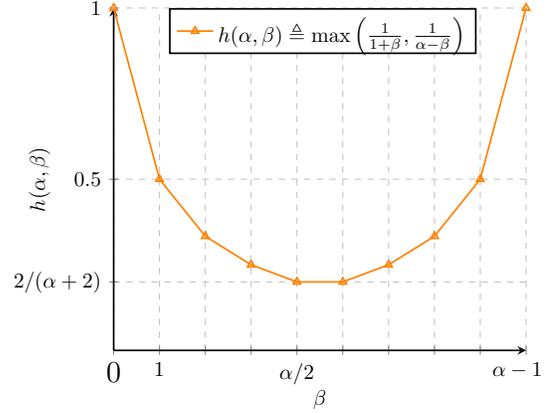
\begin{figure}[t]
 \begin{center}
 \begin{tikzpicture}[scale=0.8]
 \begin{axis}[
     axis lines = left,
     xtick = {0,1,...,9,10},
     xticklabels={\Large $0$,$1$,,,$\alpha/2$,,,,,$\!\!\!\!\!\alpha-1$},
     xlabel = $\beta$,
     ytick = {0.2,0.5,1},
     yticklabels={${2}/(\alpha+2)$,$0.5$,$1$},
     ylabel = {$h(\alpha,\beta)$},
     ymajorgrids=true,
     xmajorgrids=true,
     grid style=dashed,
     legend style={at={(0.88,1)}},
     ymin=0, ymax=1,
     line width=0.8pt,
 ]
 \addplot [
     domain=0:10, 
     samples=11, 
     color=orange,
     mark=triangle,
     ]
     {max(1/(1+x),1/(10-x))};
 \addlegendentry{$h(\alpha,\beta)\triangleq\max\left(\frac{1}{1+\beta},\frac{1}{\alpha-\beta}\right)$}
 \end{axis}
 \end{tikzpicture}
 \caption{\small Given a natural number $\alpha$, 
 plot of $h(\alpha,\beta)$ versus~$\beta$.  Given  a finite~$\alpha$, then $h(\alpha,\beta)$ is always non-zero, with minimum value $2/(\alpha+2)$, and maximum value $1$.
 }\label{fig:bounds2}
 \end{center}
 \end{figure} 
 
 \definecolor{OliveGreen}{rgb}{0,0.6,0}
 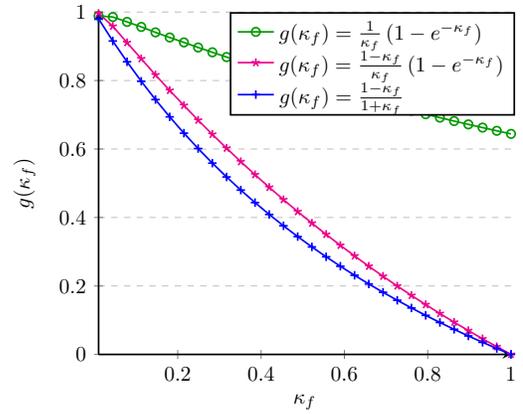
\begin{figure}[t]
 \begin{center}
 \begin{tikzpicture}[scale=0.8]
 \begin{axis}[
     axis lines = left,
     xlabel = $\kappa_f$,
     ylabel = {$g(\kappa_f)$},
     ymajorgrids=true,
     grid style=dashed,
     legend style={at={(1.01,1)}},
     ymin=0, ymax=1,
        line width=0.8pt,
 ]
       \addplot [
           domain=0.01:1, 
           samples=30, 
           color=OliveGreen,
           mark=o,
           ]
           {1/x*(1-exp(-x))+0.012};
  \addlegendentry{\!\!\!\!\!\!$g(\kappa_f)=\frac{1}{\kappa_f}\left(1-e^{-\kappa_f}\right)$}
  
 \addplot [
     domain=0.01:1, 
     samples=30, 
     color=magenta,
     mark=star,
     ]
     {1.03*(1-x)/x*(1-exp(-x))+0.00};
      \addlegendentry{$g(\kappa_f)=\frac{1-\kappa_f}{\kappa_f}\left(1-e^{-\kappa_f}\right)$}
       \addplot [
           domain=0.01:1, 
           samples=30, 
           color=blue,
           mark=+,
           ]
           {(1-x)/(1+x)};
            \addlegendentry{\hspace{-1.65cm}$g(\kappa_f)=\frac{1-\kappa_f}{1+\kappa_f}$}
 \end{axis}
 \end{tikzpicture}
 \caption{Plot of $g(\kappa_f)$ versus curvature $\kappa_f$ of a monotone submodular function $f$. By definition, the curvature $\kappa_f$ of a monotone submodular function $f$ takes values  between $0$ and $1$. $g(\kappa_f)$ increases from $0$ to $1$ as $\kappa_f$ decreases from $1$ to $0$.
 }\label{fig:bounds}
 \end{center}
 \end{figure}
 
The aforementioned classes of functions $f$ with $\kappa_f<1$ or $c_f<1$ appear in applications of control, robotics, and optimization, such as actuator and sensor placement~\cite{summers2016actuator,tzoumas2016minimal,tzoumas2016near,zhang2017kalman}, sparse approximation and feature selection~\cite{das2011spectral,khanna2017scalable}, and sparse
recovery and column subset selection~\cite{candes2006stable,boutsidis2009improved}; as a result, Problem~\ref{pr:robust_sub_max} enables critical applications such as resilient actuator placement for minimal control effort, resilient multi-robot navigation with minimal sensing and communication, and resilient experiment design; see, for example,~\cite{brent2018resilient}.
\end{myremark}
 
\paragraph{Approximation performance for low curvature}
For both monotone submodular and merely monotone functions $f$\!,  when the curvature $\kappa_f$ and the total curvature~$c_f$, respectively, tend to zero, Algorithm~\ref{alg:rob_sub_max} becomes exact,
since for $\kappa_f\rightarrow 0$ and $c_f\rightarrow 0$ the terms $\frac{1-\kappa_f}{\kappa_f}(1-e^{-\kappa_f})$, $\frac{1-\kappa_f}{1+\kappa_f}$, and $(1-c_f)^3$ in ineqs.~\eqref{ineq:bound_sub_uniform}-\eqref{ineq:bound_non_sub} respectively, tend to $1$.
Overall, Algorithm~\ref{alg:rob_sub_max}'s curvature-dependent
approximation bounds make a first step towards separating
the classes of monotone submodular and merely monotone  functions into
functions for which Problem~\ref{pr:robust_sub_max}
can be approximated well (low curvature functions), and functions for which it cannot \mbox{(high curvature functions).}

A machine learning  problem where Algorithm~\ref{alg:rob_sub_max} guarantees an approximation performance close to $100\%$ the optimal is that of Gaussian process regression for processes with RBF kernels~\cite{krause2008near,bishop2006pattern}; this problem emerges in applications of sensor deployment and scheduling for temperature monitoring.  The reason that in this class of regression problems Algorithm~\ref{alg:rob_sub_max} performs almost optimally is that the involved objective function is the entropy of the selected sensor measurements, which for Gaussian processes with RBF kernels has curvature value close to zero~\cite[Theorem~5]{sharma2015greedy}. 

\paragraph{Approximation performance for no attacks or failures}
Both for monotone submodular functions~$f$\!, and for merely monotone functions $f$\!, when the number of set removals is zero, ---i.e., when $\calI'=\emptyset$ in Problem~\ref{pr:robust_sub_max}, which implies $\beta=0$ in Theorem~\ref{th:alg_rob_sub_max_performance},--- Algorithm~\ref{alg:rob_sub_max}'s approximation performance is the same as that of the state-of-the-art algorithms for (non-resilient) set function maximization.  In~particular, \textit{for monotone submodular  functions}, scalable algorithms for (non-resilient) set function maximization have approximation performance at least $\frac{1}{\kappa_f}(1-e^{-\kappa_f})$ the optimal for any uniform matroid constraint~\cite[Theorem~5.4]{conforti1984curvature}, and $\frac{1}{1+\kappa_f}$ the optimal for any matroid constraint~\cite[Theorem~2.3]{conforti1984curvature}; at the same time, per Theorem~\ref{th:alg_rob_sub_max_performance},  when $\beta=0$, then Algorithm~\ref{alg:rob_sub_max} also has approximation performance at least $\frac{1}{\kappa_f}(1-e^{-\kappa_f})$  the optimal for any uniform matroid constraint, and $\frac{1}{1+\kappa_f}$ the optimal for any matroid constraint, since  for $\beta=0$ it is $h(\alpha,\beta)=1$ in ineq.~\eqref{ineq:bound_sub_uniform} and ineq.~\eqref{ineq:bound_sub}. Finally, \textit{for monotone  functions~$f$\!}, and for $\calI'=\emptyset$, Algorithm~\ref{alg:rob_sub_max} is the same as the algorithm proposed in~\cite[Section~2]{fisher1978analysis} for (non-resilient) set function maximization, whose  performance is optimal~\cite[Theorem~8.6]{sviridenko2017optimal}.

\myParagraph{Minimal running time}
Theorem~\ref{th:alg_rob_sub_max_performance} implies that Algorithm~\ref{alg:rob_sub_max}, even though it goes beyond the objective of (non-resilient) set function optimization, by accounting for attacks and failures, it has the same order of running time as state-of-the-art algorithms for (non-resilient) set function optimization. In particular, such algorithms for (non-resilient) set function optimization~\cite{nemhauser78analysis,fisher1978analysis,sviridenko2017optimal} terminate with $O(|\calV|^2)$ evaluations of the function~$f$\!, and Algorithm~\ref{alg:rob_sub_max} also terminates with $O(|\calV|^2)$ evaluations of the  function~$f$\!.

\myParagraph{Summary of theoretical results} In sum, Algorithm~\ref{alg:rob_sub_max} is the first algorithm for the problem of resilient maximization over matroid constraints (Problem~\ref{pr:robust_sub_max}), and it enjoys:
\begin{itemize}
\item \textit{system-wide resiliency}: Algorithm~\ref{alg:rob_sub_max} is valid for any number of denial-of-service attacks and failures;
\item \textit{minimal running time}:  Algorithm~\ref{alg:rob_sub_max} terminates with the same running time as state-of-the-art algorithms for (non-resilient) matroid-constrained optimization;
\item \textit{provable approximation performance}: for all monotone objective functions $f$ that are either submodular or merely non-decreasing with total curvature $c_f<1$, Algorithm~\ref{alg:rob_sub_max} ensures  a solution finitely close to the optimal.
\end{itemize}
Overall, Algorithm~\ref{alg:rob_sub_max} makes the first step to ensure the success of critical applications in control, robotics, and optimization~\cite{candes2006stable,boutsidis2009improved,summers2016actuator,tzoumas2016minimal,tzoumas2016near,zhang2017kalman,carlone2016attention,liu2018submodular,jawaid2015submodularity,clark2017toward,tokekar2014multi,cevher2011greedy,das2011spectral,elenberg2016restricted,khanna2017scalable,PEQUITO2017261,liu2017submodular,williams2017matroid,calinescu2007maximizing}, despite compromising worst-case attacks or failures, and with minimal running time. 

\section{Numerical Experiments on Control-Aware Sensor Selection}\label{sec:simulations}

In this section, we demonstrate the performance of Algorithm~\ref{alg:rob_sub_max} in numerical experiments.  In particular, we consider a control-aware sensor selection scenario, namely,  \textit{sensing-constrained robot navigation}, where the robot's localization for navigation is supported by on-board sensors to the robot, as well as by deployed sensors in the environment.\footnote{The scenario of {sensing-constrained robot navigation} with on-board sensors is introduced and motivated  in~\cite[Section~V]{tzoumas2018codesign}; see also~\cite{vitus2011sensor} for the case of autonomous robot navigation with deployed sensors in the environment.} 
Specifically, we consider an unmanned aerial vehicle (UAV) which has the objective to land but whose battery and measurement-processing power is limited. As a result, the UAV can to activate only a subset of its available sensors so to
localize itself, and to enable that way the generation of a control input
for landing. Specifically, we consider that the UAV generates its control input via an LQG controller, given the  measurements from the activated sensor set~\cite{bertsekas2005dynamic}.

In more detail, herein we present a Monte Carlo analysis of the above sensing-constrained robot navigation scenario for instances where sensor failures are present, and observe that Algorithm~\ref{alg:rob_sub_max} results to a near-optimal sensor selection; that is, the resulting navigation performance of the UAV matches the optimal in all tested instances where the optimal sensor selection could be computed via a brute-force algorithm.

\myParagraph{Simulation setup} We consider an UAV that moves in a 3D space, starting from a
randomly selected initial location. 
The objective of the UAV is to land at 
 position $[0,\;0,\;0]$ with zero velocity. 
The UAV is modelled as a double-integrator
with state $x_t = [p_t \; v_t]^\top \in \Real{6}$ at each time $t=1,2,\ldots$  
 ($p_t$ is the 3D position of the UAV, and $v_t$ is its velocity), and can control its own acceleration 
$u_t \in \Real{3}$; the process noise is chosen as $W_t = \eye_6$. 
The UAV may support its localization by utilizing $2$ on-board sensors and $12$ deployed sensors on the ground. The on-board sensors are one GPS receiver, measuring the 
UAV position  $p_t$ with a covariance~$2 \cdot\eye_3$, 
and one altimeter, measuring only the last component of $p_t$ (altitude) with standard deviation $0.5\rm{m}$. The ground sensors vary with each Monte Carlo run, and are generated randomly; we consider them to provide linear measurements of the UAV's state.
Among the aforementioned $14$ available sensors to the UAV, we assume that the UAV can use only $\alpha$ of them.

In particular, the UAV chooses the $\alpha$ sensors to activate so to minimize an LQG cost of the form:
\begin{equation}\label{eq:lqg_cost}
\sum_{t=1}^{T}[x_t^\top Qx_t+u_t^\top Ru_t],
\end{equation}
per the problem formulation in~\cite[Section~II]{tzoumas2018codesign}, where the cost matrix $Q$ penalizes the deviation of the state vector from the zero state (since the UAV's objective is to land at position $[0,\;0,\;0]$ with zero velocity), and the cost matrix $R$ penalizes the control input vector;
specifically, in the simulation setup herein we consider $Q = \diag{[1e^{-3},\; 1e^{-3},\;10,\; 1e^{-3},\; 1e^{-3},\; 10]}$ 
and $R = \eye_3$.  Note that the structure of $Q$ reflects the fact that during landing 
we are particularly interested in controlling the vertical direction and the vertical velocity 
(entries with larger weight in $Q$), while we are less interested in controlling accurately the 
horizontal position and velocity (assuming a sufficiently large landing site). In~\cite[Section~III]{tzoumas2018codesign} it is proven that the UAV selects an optimal sensor set $\calS$, and enables the generation of an optimal LQG control input with cost matrices $Q$ and $R$, if it selects $\calS$ by minimizing an objective function of the form:
\begin{equation}\label{eq:opt_sensors}
\sum_{t=1}^{T}\text{trace}[M_t\Sigma_{t|t}(\calS)],
\end{equation}
where $M_t$ is a positive semi-definite matrix that depends on the LQG cost matrices $Q$ and $R$, as well as, on the UAV's system dynamics; and $\Sigma_{t|t}(\calS)$ is the error covariance of the Kalman filter given the sensor set selection $\calS$.

\myParagraph{Compared algorithms}
We compare four algorithms; all algorithms
only differ in how they select the sensors used.
The~first algorithm is the optimal sensor selection algorithm, denoted as \toptimal, which 
attains the minimum of the cost function in eq.~\eqref{eq:opt_sensors}; this brute-force approach is viable since the number of available sensors is small.
The second approach is a random sensor selection, denoted as {\tt random$^*$}\!.
The third approach, denoted as \tlogdet, selects sensors to greedily minimize the cost function in eq.~\eqref{eq:opt_sensors}, \textit{ignoring the possibility of sensor failures}, per the problem formulation in eq.~\eqref{eq:non_res}.
The fourth approach uses Algorithm~\ref{alg:rob_sub_max} to solve the resilient re-formulation of eq.~\eqref{eq:opt_sensors} per Problem~\ref{pr:robust_sub_max}, and is denoted as \tslqg.  From each of the selected sensor sets, by each of the above four algorithms respectively, we consider an optimal sensor removal, which we compute via brute-force.

\myParagraph{Results}  We next present our simulation results averaged over 20 Monte Carlo runs of the above simulation setup, where we vary the number~$\alpha$ of sensor selections from $2$ up to $12$ with step $1$, and the number~$\beta$ of sensors failures from $1$ to $10$  with step $3$, and where we~randomize the 
sensor matrices of the $12$ ground sensors.  In particular, the results of our numerical analysis are reported in Fig.~\ref{fig:formationControlStats}.  
In more detail, Fig.~\ref{fig:formationControlStats} shows the attained \LQG cost for all the combinations of $\alpha$ and $\beta$ values where $\beta\leq \alpha$ (for $\beta>\alpha$ the LQG cost is considered $+\infty$, since $\beta>\alpha$ implies that all $\alpha$ selected sensors fail).  The following observations from Fig.~\ref{fig:formationControlStats} are due: 
\begin{itemize}
\item \textit{Near-optimality of the \tslqg algorithm (Algorithm~\ref{alg:rob_sub_max})}: Algorithm~\ref{alg:rob_sub_max} ---blue colour in Fig.~\ref{fig:formationControlStats}--- performs close to the optimal algorithm \toptimal ---green colour in Fig.~\ref{fig:formationControlStats}. In particular, across all but two scenarios in Fig.~\ref{fig:formationControlStats}, Algorithm~\ref{alg:rob_sub_max} achieves an approximation performance at least 97\% the optimal; and in the remaining two scenarios (see Fig.~\ref{fig:formationControlStats}-(a) for $\alpha$ equal to $3$ or $4$, and $\beta$ equal to~$1$), Algorithm~\ref{alg:rob_sub_max} achieves an approximation performance at least 90\% the optimal.
\item \textit{Performance of the \tlogdet algorithm}: The \tlogdet algorithm ---red colour in Fig.~\ref{fig:formationControlStats}--- performs poorly as the number~$\beta$ of sensor failures increases, which is expected given that the \tlogdet algorithm minimizes the cost function in eq.~\eqref{eq:opt_sensors} {ignoring the possibility of sensor failures}.  Notably, for some of the cases the \tlogdet performs worse or equally poor as the {\tt random$^*$}: for example, see Fig.~\ref{fig:formationControlStats}-(c) for $\alpha\geq 9$, and Fig.~\ref{fig:formationControlStats}-(d).
\item \textit{Performance of the {\tt random$^*$} algorithm}: Expectedly, the performance of also the {\tt random$^*$} algorithm ---black colour in Fig.~\ref{fig:formationControlStats}--- is poor across all scenarios in Fig.~\ref{fig:formationControlStats}.
\end{itemize}

 \definecolor{OliveGreen}{rgb}{0,0.6,0}
\newcommand{\myhspace}{\hspace{2mm}}
\newcommand{\mpw}{4.5cm}
\begin{figure*}[t]
\begin{center}
\begin{minipage}{\textwidth}
\centering
\hspace{-5mm}
\begin{tabular}{cc}%
\begin{minipage}{\mpw}%
\centering
\begin{tikzpicture}[scale=0.7]
\begin{axis}[
    axis lines = left,
    ymin=2300,
    xlabel = {\large$\alpha$},
    xticklabels={$2$,$3$,$4$,$5$,$6$,$7$,$8$,$9$,$10$,$11$,$12$},
    xtick = {2,...,12},
    ylabel = {\large\text{LQG cost per eq.~\eqref{eq:lqg_cost}}},
    legend pos=south east,
    ymajorgrids=true,
    grid style=dashed,
    line width=1.1pt,
    legend style={at={(1,.7)}},
]

\addplot[
    color=black,
    mark=cross,
    style={solid},mark=star,
    ]
    coordinates { 
    (2,4.429e+04)(3,9763)(4,2.28e+04
    )(5,6033)(6,7847)(7,6630)(8,5142)(9,3311)(10,4990)(11,4163)(12,3396)
    };
\addlegendentry{{\tt random$^*$}}
\addplot[
    color=red,
    mark=cross, 
    style={densely  dotted}, mark=otimes*
    ]
    coordinates {  
     (2,1.3e+04)(3,6750)(4,5700
     )(5,5401)(6,4506)(7,4876)(8,3793)(9,2469)(10,3128)(11,4063)(12,3382)
    };
\addlegendentry{\tlogdet}
\addplot[
    color=blue,
    mark=cross,
    style={solid},mark=square
    ]
    coordinates {  
     (2,8874)(3,6909)(4,6038
     )(5,4964)(6,4298)(7,4876)(8,3793)(9,2469)(10,3128)(11,4063)(12,3382)
    };
\addlegendentry{\tslqg}
\addplot[
    color=OliveGreen,
    style=solid,
    ]
    coordinates {
     (2,8874)(3,5983)(4,4907
     )(5,4964)(6,4298)(7,4876)(8,3793)(9,2469)(10,3128)(11,4063)(12,3382)
    };
\addlegendentry{\toptimal}
\end{axis}
\end{tikzpicture} \\
(a) \hspace{1.9cm}$\beta=1$
\end{minipage}
& \hspace{20mm}
\begin{minipage}{\mpw}%
\centering%
\begin{tikzpicture}[scale=0.7]
\begin{axis}[
    axis lines = left,
    ymin=3000,
    xlabel = {\large$\alpha$},
        xticklabels={$2$,$3$,$4$,$5$,$6$,$7$,$8$,$9$,$10$,$11$,$12$},
        xtick = {2,...,12},
    ylabel = {\large\text{LQG cost per eq.~\eqref{eq:lqg_cost}}},
    legend pos=south east,
    ymajorgrids=true,
    grid style=dashed,
    line width=1.1pt,
         legend style={at={(1,.7)}},
]

\addplot[
    color=black,
    mark=cross,
    style={solid},mark=star,
    ]
    coordinates { 
    (5,3.008e+05)(6,5.906e+04)(7,4.099e+04
    )(8,1.974e+04)(9,1.217e+04)(10,1.392e+04)(11,1.414e+04)(12,1.258e+04)
    };
\addlegendentry{{\tt random$^*$}}
\addplot[
    color=red,
    mark=cross, 
    style={densely  dotted}, mark=otimes*
    ]
    coordinates {  
    (5,1.066e+05)(6,3.244e+04)(7,9947
    )(8,9485)(9,1.063e+04)(10,1.251e+04)(11,1.181e+04)(12,1.021e+04)
    };
\addlegendentry{\tlogdet}
\addplot[
    color=blue,
    mark=cross,
    style={solid},mark=square
    ]
    coordinates {  
    (5,4.285e+04)(6,1.671e+04)(7,9354
    )(8,7411)(9,1.063e+04)(10,1.251e+04)(11,1.072e+04)(12,9991)
    };
\addlegendentry{\tslqg}
\addplot[
    color=OliveGreen,
    style=solid,
    ]
    coordinates {
    (5,4.285e+04)(6,1.671e+04)(7,9354
    )(8,7411)(9,1.063e+04)(10,1.251e+04)(11,1.072e+04)(12,9991)
    };
\addlegendentry{\toptimal}

\end{axis}
\end{tikzpicture} \\
(b) \hspace{1.9cm}$\beta=4$ 
\end{minipage}\\
\myhspace
\\
\begin{minipage}{\mpw}%
\centering
\begin{tikzpicture}[scale=0.7]
\begin{axis}[
    axis lines = left,
    ymin=2300,
    xlabel = {\large$\alpha$},
        xticklabels={$2$,$3$,$4$,$5$,$6$,$7$,$8$,$9$,$10$,$11$,$12$},
        xtick = {2,...,12},
    ylabel = {\large\text{LQG cost per eq.~\eqref{eq:lqg_cost}}},
    legend pos=south east,
    ymajorgrids=true,
    grid style=dashed,
    line width=1.1pt,
         legend style={at={(1,.7)}},
]

\addplot[
    color=black,
    mark=cross,
    style={solid},mark=star,
    ]
    coordinates { 
 (8,3.554e+05)(9,4.797e+04)(10,4.166e+04)(11,3.007e+04)(12,2.578e+04)
    };
\addlegendentry{{\tt random$^*$}}
\addplot[
    color=red,
    mark=cross, 
    style={densely  dotted}, mark=otimes*
    ]
    coordinates {  
 (8,1.853e+05)(9,5.837e+04)(10,4.465e+04)(11,3.176e+04)(12,2.326e+04)
    };
\addlegendentry{\tlogdet}
\addplot[
    color=blue,
    mark=cross,
    style={solid},mark=square
    ]
    coordinates {  
 (8,8.034e+04)(9,3.527e+04)(10,3.545e+04)(11,2.033e+04)(12,2.326e+04)
    };
\addlegendentry{\tslqg}
\addplot[
    color=OliveGreen,
    style=solid,
    ]
    coordinates {
 (8,8.034e+04)(9,3.488e+04)(10,3.545e+04)(11,2.033e+04)(12,2.326e+04)
    };
\addlegendentry{\toptimal}

\end{axis}
\end{tikzpicture} \\
(c) \hspace{1.9cm}$\beta=7$
\end{minipage}
&  \hspace{20mm}
\begin{minipage}{\mpw}%
\centering%
\begin{tikzpicture}[scale=0.7]
\begin{axis}[
    axis lines = left,
    ymin=2300,
    xlabel = {\large$\alpha$},
        xticklabels={$2$,$3$,$4$,$5$,$6$,$7$,$8$,$9$,$10$,$11$,$12$},
        xtick = {2,...,12},
    ylabel = {\large\text{LQG cost per eq.~\eqref{eq:lqg_cost}}},
    legend pos=south east,
    ymajorgrids=true,
    grid style=dashed,
    line width=1.1pt,
         legend style={at={(1,.7)}},
]

\addplot[
    color=black,
    mark=cross,
    style={solid},mark=star,
    ]
    coordinates { 
(11,7.352e+05)(12,6.833e+04)
    };
\addlegendentry{{\tt random$^*$}}
\addplot[
    color=red,
    mark=cross, 
    style={densely  dotted}, mark=otimes*
    ]
    coordinates {  
(11,7.352e+05)(12,5.36e+04)
    };
\addlegendentry{\tlogdet}
\addplot[
    color=blue,
    mark=cross,
    style={solid},mark=square
    ]
    coordinates {  
(11,1.624e+05)(12,5.36e+04)
    };
\addlegendentry{\tslqg}
\addplot[
    color=OliveGreen,
    style=solid,
    ]
    coordinates {
(11,1.624e+05)(12,5.36e+04)
    };
\addlegendentry{\toptimal}

\end{axis}
\end{tikzpicture} \\
(d) \hspace{2.05cm}$\beta=10$ 
\end{minipage}
\end{tabular}
\end{minipage}%
\vspace{-1mm}
\caption{\label{fig:formationControlStats}\small
\LQG cost for increasing number of sensor selections $\alpha$ (from $2$ up to $12$ with step $1$), and for  $4$ values of $\beta$ (number of sensor failures among the $\alpha$ selected sensors); in particular, the value of $\beta$ varies across the sub-figures as follows: $\beta=1$ in sub-figure (a); $\beta=4$ in sub-figure (b); $\beta=7$ in sub-figure (c); and $\beta=10$ in sub-figure (d).
}\vspace{-5mm}
\end{center}
\end{figure*}
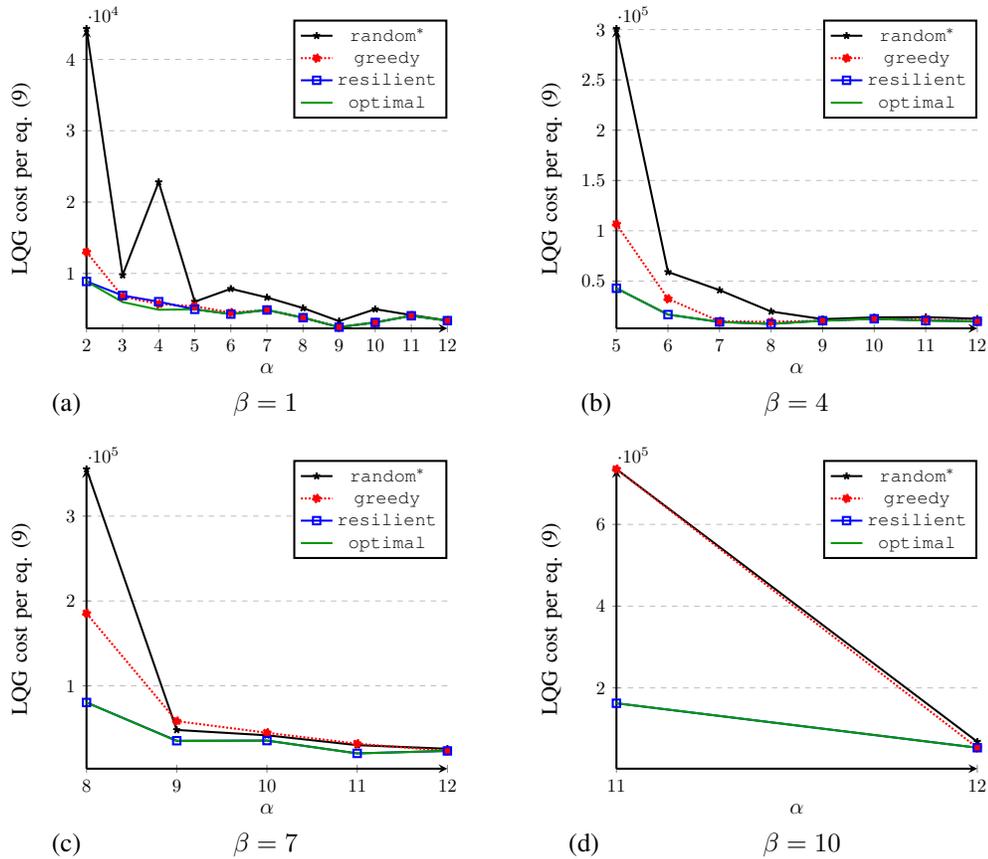

Overall, in the above numerical experiments, Algorithm~\ref{alg:rob_sub_max} demonstrates a close-to-optimal approximation performance, and the necessity for a resilient re-formulation of the optimization problem in eq.~\eqref{eq:non_res}, e.g., per Problem~\ref{pr:robust_sub_max}, is exemplified.

\section{Concluding remarks \& Future work} \label{sec:con}

We made the first step to ensure the success of critical missions in control, robotics, and optimization that involve the design of systems subject to complex optimization constraints of heterogeneity and interdependency ---called matroid constraints--- against worst-case denial-of-service attacks or failures.  In particular,  we provided the first algorithm for Problem~\ref{pr:robust_sub_max}, which, with minimal running time, guarantees
a close-to-optimal performance against system-wide attacks and failures.  To quantify the algorithm's approximation performance, we exploited a notion of curvature for monotone (not necessarily submodular) set functions, and contributed a first step towards characterizing the curvature's effect on the approximability of resilient \textit{matroid-constrained} maximization.  Our curvature-dependent characterizations complement the current knowledge on the curvature's effect on the approximability of simpler problems, such as of {non-matroid-constrained} resilient maximization~\cite{tzoumas2017resilient,bogunovic2018robust,tzoumas2018resilientSequential}, and of {non-resilient} maximization~\cite{conforti1984curvature,iyer2013curvature,bian2017guarantees}. Finally, we supported our theoretical analyses with numerical experiments.

This paper opens several avenues for future research, both in theory and in applications.
Future work in theory includes the extension of our results to sequential (multi-step) maximization, per the recent developments in~\cite{tzoumas2018resilientSequential}, to enable applications of sensor scheduling and of path planning in online optimization that \textit{adapts} against persistent attacks and failures~\cite{tokekar2014multi,tzoumas2016scheduling}.
Future work in applications includes the experimental testing of the proposed algorithm in applications of motion-planning for multi-target covering with mobile vehicles~\cite{tokekar2014multi}, 
to~enable resiliency in {critical scenarios of surveillance.

\section{Acknowledgements}

We thank Luca Carlone and Konstantinos Gatsis for inspiring discussions, and Luca Carlone for sharing code that enabled the numerical experiments in this paper.

\bibliographystyle{IEEEtran}
\bibliography{references,security_references}

\appendices

\section{Notation}\label{app:notation}

In the appendices below we use the following notation: given a finite ground set $\calV$, and a set function $f:2^\mathcal{V}\mapsto \mathbb{R}$, then, for any sets $\mathcal{X}\subseteq \mathcal{V}$ and $\mathcal{X}'\subseteq \mathcal{V}$: 
\begin{equation*}\label{notation:marginal}
f(\mathcal{X}|\mathcal{X}')\triangleq f(\mathcal{X}\cup\mathcal{X}')-f(\mathcal{X}').
\end{equation*}
Moreover,  let the set $\mathcal{A}^\star$ denote an (optimal) solution to Problem~\ref{pr:robust_sub_max}; formally:
\begin{equation*}
\mathcal{A}^\star\in \arg\underset{\mathcal{A}\subseteq \calV, \calA\in \calI}{\max} \; \; \underset{\;\mathcal{B}\subseteq \calA, \calB\in \calI'(\calA)}{\min}   \; \;\; f(\mathcal{A}\setminus \mathcal{B}).
\end{equation*}

\section{Preliminary lemmas}\label{app:prelim}

We list lemmas that support the proof of Theorem~\ref{th:alg_rob_sub_max_performance}.\footnote{The proof of Lemmas~\ref{lem:D3}-\ref{lem:subratio} and of Corollary~\ref{cor:ineq_from_lemmata} is also found in~\cite{tzoumas2017resilient} and~\cite{tzoumas2018resilientSequential}.}

\begin{mylemma}\label{lem:D3}
Consider any finite ground set $\mathcal{V}$, a non-decreasing submodular function $f:2^\mathcal{V}\mapsto \mathbb{R}$, and non-empty sets $\calY, \calP \subseteq \calV$ such that for all elements $y \in \calY$, and all elements $p \in \calP$, it is $f(y)\geq f(p)$.  Then:
\belowdisplayskip=-12pt\begin{equation*}
f(\calP|\calY)\leq |\calP|f(\calY).
\end{equation*}
\end{mylemma}
\paragraph*{Proof of Lemma~\ref{lem:D3}} Consider any element $y \in \calY$; then: 
\begin{align}
f(\calP|\calY)&= f(\calP\cup\calY)-f(\calY)\label{aux1:1}\\
&\leq f(\calP)+f(\calY)-f(\calY)\label{aux1:2}\\
&= f(\calP)\nonumber\\
&\leq \sum_{p\in\calP}f(p)\label{aux1:4}\\
&\leq |\calP| \max_{p\in\calP} f(p)\nonumber\\
&\leq |\calP|  f(y)\label{aux1:7}\\
&\leq |\calP| f(\calY),\label{aux1:5}
\end{align}
where eqs.~\eqref{aux1:1}-\eqref{aux1:5} hold for the following reasons: eq.~\eqref{aux1:1} holds since for any sets $\mathcal{X}\subseteq \mathcal{V}$ and $\mathcal{Y}\subseteq \mathcal{V}$, it is $f(\mathcal{X}|\mathcal{Y})=f(\mathcal{X}\cup \mathcal{Y})-f(\mathcal{Y})$; ineq.~\eqref{aux1:2} holds since $f$ is submodular and, as a result, the submodularity Definition~\ref{def:sub} implies that for any set $\calA\subseteq\calV$ and $\calA'\subseteq\calV$, it is $f(\calA\cup \calA')\leq f(\calA)+f(\calA')$~\cite[Proposition 2.1]{nemhauser78analysis}; ineq.~\eqref{aux1:4} holds for the same reason as ineq.~\eqref{aux1:2}; ineq.~\eqref{aux1:7} holds since for all elements $y \in \calY$, and for all elements $p \in \calP$, it is $f(y)\geq f(p)$; finally, ineq.~\eqref{aux1:5} holds since $f$ is monotone, and since $y\in\calY$.
\hfill $\blacksquare$

\begin{mylemma}\label{lem:non_total_curvature}
Consider a finite ground set $\mathcal{V}$, and a non-decreasing submodular set function $f:2^\mathcal{V}\mapsto \mathbb{R}$ such that $f$ is non-negative and $f(\emptyset)=0$. Then, for any $\mathcal{A}\subseteq \mathcal{V}$, it~holds:
\begin{equation*}
f(\mathcal{A})\geq (1-\kappa_f)\sum_{a \in \mathcal{A}}f(a).
\end{equation*}
\end{mylemma}
\paragraph*{Proof of Lemma~\ref{lem:non_total_curvature}} Let $\mathcal{A}=\{a_1,a_2,\ldots, a_{|{\cal A}|}\}$. We prove Lemma~\ref{lem:non_total_curvature} by proving the following two inequalities: 
\begin{align}
f(\mathcal{A})&\geq \sum_{i=1}^{|{\cal A}|} f(a_i|\mathcal{V}\setminus\{a_i\}),\label{ineq5:aux_5}\\
\sum_{i=1}^{|{\cal A}|} f(a_i|\mathcal{V}\setminus\{a_i\})&\geq (1-\kappa_f)\sum_{i=1}^{|{\cal A}|} f(a_i)\label{ineq5:aux_6}. 
\end{align} 

We begin with the proof of ineq.~\eqref{ineq5:aux_5}: 
\begin{align}
f(\mathcal{A})&=f(\mathcal{A}|\emptyset)\label{ineq5:aux_9}\\
&\geq f(\mathcal{A}|\mathcal{V}\setminus \mathcal{A})\label{ineq5:aux_10}\\
&= \sum_{i=1}^{|{\cal A}|}f(a_i|\mathcal{V}\setminus\{a_i,a_{i+1},\ldots,a_{|{\cal A}|}\})\label{ineq5:aux_11}\\
&\geq \sum_{i=1}^{|{\cal A}|}f(a_i|\mathcal{V}\setminus\{a_i\}),\label{ineq5:aux_12}
\end{align}
where ineqs.~\eqref{ineq5:aux_10}-\eqref{ineq5:aux_12} hold for the following reasons: ineq.~\eqref{ineq5:aux_10} is implied by eq.~\eqref{ineq5:aux_9} because $f$ is submodular and $\emptyset\subseteq \mathcal{V}\setminus \mathcal{A}$; eq.~\eqref{ineq5:aux_11} holds since for any sets $\mathcal{X}\subseteq \mathcal{V}$ and $\mathcal{Y}\subseteq \mathcal{V}$ it is $f(\mathcal{X}|\mathcal{Y})=f(\mathcal{X}\cup \mathcal{Y})-f(\mathcal{Y})$, and since $\{a_1,a_2,\ldots, a_{|{\cal A}|}\}$ denotes the set $\mathcal{A}$; and ineq.~\eqref{ineq5:aux_12} holds since $f$ is submodular, and since $\mathcal{V}\setminus\{a_i,a_{i+1},\ldots,a_{\mu}\} \subseteq \mathcal{V}\setminus\{a_i\}$.  These observations complete the proof of ineq.~\eqref{ineq5:aux_5}.

We now prove ineq.~\eqref{ineq5:aux_6} using the Definition~\ref{def:curvature} of $\kappa_f$, as follows: since $\kappa_f=1-\min_{v\in \mathcal{V}}\frac{f(v|\mathcal{V}\setminus\{v\})}{f(v)}$, it is implied that for all elements $v\in \mathcal{V}$ it is $ f(v|\mathcal{V}\setminus\{v\})\geq (1-\kappa_f)f(v)$.  Therefore, by adding the latter inequality across all elements $a \in \calA$, we complete the proof of ineq.~\eqref{ineq5:aux_6}.
\hfill $\blacksquare$

\begin{mylemma}\label{lem:curvature2}
Consider a finite ground set $\mathcal{V}$, and a monotone set function $f:2^\mathcal{V}\mapsto \mathbb{R}$ such that $f$ is non-negative and $f(\emptyset)=0$. Then, for any sets $\mathcal{A}\subseteq \mathcal{V}$ and $\mathcal{B}\subseteq \mathcal{V}$ such that $\calA \cap \calB=\emptyset$, it holds:
\begin{equation*}
f(\mathcal{A}\cup \mathcal{B})\geq (1-c_f)\left(f(\mathcal{A})+f(\mathcal{B})\right).
\end{equation*}
\end{mylemma}
\paragraph*{Proof of Lemma~\ref{lem:curvature2}}
Let $\mathcal{B}=\{b_1, b_2, \ldots, b_{|\mathcal{B}|}\}$. Then, 
\begin{equation}
f(\mathcal{A}\cup \mathcal{B})=f(\mathcal{A})+\sum_{i=1}^{|\mathcal{B}|}f(b_i|\mathcal{A}\cup \{b_1, b_2, \ldots, b_{i-1}\}). \label{eq1:lemma_curvature2}
\end{equation} 
The definition of total curvature in Definition~\ref{def:total_curvature} implies:
\begin{align}
&\!\!\!f(b_i|\mathcal{A}\cup \{b_1, b_2, \ldots, b_{i-1}\})\geq\nonumber\\
& (1-c_f)f(b_i|\{b_1, b_2, \ldots, b_{i-1}\}). \label{eq2:lemma_curvature2}
\end{align} 
The proof is completed by substituting ineq.~\eqref{eq2:lemma_curvature2} in eq.~\eqref{eq1:lemma_curvature2} and then by taking into account that it holds $f(\mathcal{A})\geq (1-c_f)f(\mathcal{A})$, since $0\leq c_f\leq 1$.
\hfill $\blacksquare$

\begin{mylemma}\label{lem:curvature}
Consider a finite ground set $\mathcal{V}$m and a non-decreasing set function $f:2^\mathcal{V}\mapsto \mathbb{R}$ such that $f$ is non-negative and $f(\emptyset)=0$. Then, for any set $\mathcal{A}\subseteq \mathcal{V}$ and any set $\mathcal{B}\subseteq \mathcal{V}$ such that $\calA \cap \calB=\emptyset$, it holds:
\begin{equation*}
f(\mathcal{A}\cup \mathcal{B})\geq (1-c_f)\left(f(\mathcal{A})+\sum_{b \in \mathcal{B}}f(b)\right).
\end{equation*}
\end{mylemma}
\paragraph*{Proof of Lemma~\ref{lem:curvature}}
Let $\mathcal{B}=\{b_1, b_2, \ldots, b_{|\mathcal{B}|}\}$. Then, 
\begin{equation}
f(\mathcal{A}\cup \mathcal{B})=f(\mathcal{A})+\sum_{i=1}^{|\mathcal{B}|}f(b_i|\mathcal{A}\cup \{b_1, b_2, \ldots, b_{i-1}\}). \label{eq1:lemma_curvature}
\end{equation} 
In addition, Definition~\ref{def:total_curvature} of total curvature implies:
\begin{align}
f(b_i|\mathcal{A}\cup \{b_1, b_2, \ldots, b_{i-1}\})&\geq (1-c_f)f(b_i|\emptyset)\nonumber\\
&=(1-c_f)f(b_i), \label{eq2:lemma_curvature}
\end{align} 
where the latter equation holds since $f(\emptyset)=0$.
The proof is completed by substituting~\eqref{eq2:lemma_curvature} in~\eqref{eq1:lemma_curvature} and then taking into account that $f(\mathcal{A})\geq (1-c_f)f(\mathcal{A})$ since $0\leq c_f\leq 1$. \hfill $\blacksquare$

\begin{mylemma}\label{lem:subratio}
Consider a finite ground set $\mathcal{V}$ and a non-decreasing set function $f:2^\mathcal{V}\mapsto \mathbb{R}$ such that $f$ is non-negative and $f(\emptyset)=0$. Then, for any set $\mathcal{A}\subseteq \mathcal{V}$ and any set  $\mathcal{B}\subseteq \mathcal{V}$ such that $\mathcal{A}\setminus\mathcal{B}\neq \emptyset$, it holds:
\begin{equation*}
f(\mathcal{A})+(1-c_f) f(\mathcal{B})\geq (1-c_f) f(\mathcal{A}\cup \mathcal{B})+f(\mathcal{A}\cap \mathcal{B}).
\end{equation*}
\end{mylemma}
\paragraph*{Proof of Lemma~\ref{lem:subratio}}
Let $\mathcal{A}\setminus\mathcal{B}=\{i_1,i_2,\ldots, i_r\}$, where $r=|\mathcal{A}-\mathcal{B}|$. From Definition~\ref{def:total_curvature} of total curvature $c_f$, for any $i=1,2, \ldots, r$\!, it is  $f(i_j|\mathcal{A} \cap \mathcal{B} \cup \{i_1, i_2, \ldots, i_{j-1}\})\geq (1-c_f) f(i_j|\mathcal{B} \cup \{i_1, i_2, \ldots, i_{j-1}\})$. Summing these $r$ inequalities:
$$f(\mathcal{A})-f(\mathcal{A}\cap \mathcal{B})\geq (1-c_f) \left(f(\mathcal{A}\cup \mathcal{B})-f(\mathcal{B})\right),$$
which implies the lemma. \hfill $\blacksquare$

\begin{mycorollary}\label{cor:ineq_from_lemmata}
Consider a finite ground set $\mathcal{V}$m and a non-decreasing set function $f:2^\mathcal{V}\mapsto \mathbb{R}$ such that $f$ is non-negative and $f(\emptyset)=0$. Then, for any set $\mathcal{A}\subseteq \mathcal{V}$ and any set $\mathcal{B}\subseteq \mathcal{V}$ such that $\mathcal{A}\cap\mathcal{B}=\emptyset$, it holds:
\begin{equation*}
f(\mathcal{A})+\sum_{b \in \mathcal{B}}f(b) \geq (1-c_f)  f(\mathcal{A}\cup \mathcal{B}).
\end{equation*}
\end{mycorollary}
\paragraph*{Proof of Corollary~\ref{cor:ineq_from_lemmata}}
 Let $\mathcal{B}=\{b_1,b_2,\ldots,b_{|\mathcal{B}|}\}$. 
\begin{align}
f(\mathcal{A})+\sum_{i=1}^{|\mathcal{B}|}f(b_i) &\geq (1-c_f) f(\mathcal{A})+\sum_{i=1}^{|\mathcal{B}|}f(b_i))\label{ineq:cor_aux1} \\
& \geq (1-c_f) f(\mathcal{A}\cup \{b_1\})+\sum_{i=2}^{|\mathcal{B}|}f(b_i)\nonumber\\
& \geq (1-c_f) f(\mathcal{A}\cup \{b_1,b_2\})+\sum_{i=3}^{|\mathcal{B}|}f(b_i)\nonumber\\
& \;\;\vdots \nonumber\\
& \geq (1-c_f) f(\mathcal{A}\cup \mathcal{B}),\nonumber
\end{align}
where~\eqref{ineq:cor_aux1} holds since $0\leq c_f\leq 1$, and the rest due to Lemma~\ref{lem:subratio}, since $\mathcal{A}\cap\mathcal{B}=\emptyset$ implies $\mathcal{A}\setminus \{b_1\}\neq \emptyset$, $\mathcal{A}\cup \{b_1\}\setminus \{b_2\}\neq \emptyset$, $\ldots$, $\mathcal{A}\cup \{b_1,b_2,\ldots, b_{|\mathcal{B}|-1}\}\setminus \{b_{|\mathcal{B}|}\}\neq \emptyset$. 

\hfill $\blacksquare$

\begin{mylemma}\label{lem:dominance}
Recall the notation in Algorithm~\ref{alg:rob_sub_max}, and consider the sets $\calA_1$ and $\calA_2$ constructed by Algorithm~\ref{alg:rob_sub_max}'s lines~\ref{line:begin_while_1}-\ref{line:end_while_1} and lines~\ref{line:begin_while_2}-\ref{line:end_while_2}, respectively. Then, for all elements $v\in \calA_1$ and all elements $v'\in \calA_2$, it holds $f(v)\geq f(v')$.
\end{mylemma}
\paragraph*{Proof of Lemma~\ref{lem:dominance}}
Let $v_1,\ldots,v_{|\calA_1|}$ be the elements in~$\calA_{1}$, ---i.e., $\calA_1\equiv \{v_1,\ldots,v_{|\calA_1|}\}$,--- and be such that for each $i=1,\ldots,|\calA_1|$ the element~$v_i$ is the $i$-th element added in $\calA_1$ per Algorithm~\ref{alg:rob_sub_max}'s lines~\ref{line:begin_while_1}-\ref{line:end_while_1}; similarly, let $v_1',\ldots,v_{|\calA_2|}'$ be the elements in $\calA_{2}$, ---i.e., $\calA_2\equiv \{v_1',\ldots,v_{|\calA_2|}'\}$,--- and be such that for each $i=1,\ldots,|\calA_2|$ the element $v_i'$ is the $i$-th element added in $\calA_2$ per Algorithm~\ref{alg:rob_sub_max}'s lines~\ref{line:begin_while_2}-\ref{line:end_while_2}.

We prove Lemma~\ref{lem:dominance} by the method of contradiction; specifically, we focus on the case where $(\calV,\calI')$ is a uniform matroid; the case where $(\calV,\calI')$ is a partition matroid with the same partition as $(\calV,\calI)$ follows the same steps by focusing at each partition separately.  In~particular, assume that there exists an index $i\in\{1,\ldots,|\calA_1|\}$ and an $j\in\{1,\ldots,|\calA_2|\}$ such that $f(v'_j)>f(v_i)$, and, in particular, assume that $i,j$ are the smallest indexes such that $f(v'_j)>f(v_i)$.  Since Algorithm~\ref{alg:rob_sub_max} constructs $\calA_1$ and~$\calA_2$ such that $\calA_1\cup\calA_2\in \calI$, and since it also is that $(\calV,\calI)$ is a matroid and $\{v_1,\ldots,v_{i-1}, v'_{j}\}\subseteq \calA_1\cup\calA_2$, we have that $\{v_1,\ldots,v_{i-1}, v'_{j}\}\in \calI$. In addition, we have that $\{v_1,\ldots,v_{i-1}, v'_{j}\}\in \calI'$\!\!, since $\calI'$ is either a uniform or a partition matroid and, as a result, if $\{v_1,\ldots,v_{i-1}, v_{i}\}\in \calI'$ then it also is $\{v_1,\ldots,v_{i-1}, v\}\in \calI'$ for any $v\in\calV\setminus\{v_1,\ldots,v_{i-1}\}$.  Overall, $\{v_1,\ldots,v_{i-1}, v'_{j}\}\in \calI,\calI'$\!\!. Now, consider the ``while loop'' in Algorithm~\ref{alg:rob_sub_max}'s lines~\ref{line:begin_while_1}-\ref{line:end_while_1} at the beginning of its $i$-th iteration, that is, when Algorithm~\ref{alg:rob_sub_max} has chosen only the elements $\{v_1,\ldots,v_{i-1}
\}$ among the elements in $\calA_1$.  Then, per Algorithm~\ref{alg:rob_sub_max}'s lines~\ref{line:select_element_bait}-\ref{line:build_of_bait}, the next element $v$ that is added in $\{v_1,\ldots,v_{i-1}
\}$ is the one that achieves the highest value of $f(v')$ among all elements in $v'\in\calV\setminus\{v_1,\ldots,v_{i-1}\}$, and which satisfies $\{v_1,\ldots,v_{i-1}, v\}\in \calI,\calI'$\!\!. Therefore, the next element $v$ that is added in $\{v_1,\ldots,v_{i-1}
\}$ cannot be $v_i$, since $f(v'_j)>f(v_i)$ and $\{v_1,\ldots,v_{i-1}, v'_{j}\}\in \calI,\calI'$.
\hfill $\blacksquare$

\begin{mylemma}\label{lem:its_a_matroid}
Consider a matroid $(\calV,\calI)$, and a set $\calY\subseteq \calV$ such that $\calY \in \calI$.  Moreover, define the following collection of subsets of $\calV\setminus \calY$: $\calI'\triangleq\{\calX: \calX \subseteq \calV\setminus\calY, \calX\cup \calY \in \calI\}$.  Then, $(\calV\setminus\calY, \calI')$ is a matroid.
\end{mylemma}
\paragraph*{Proof of Lemma~\ref{lem:its_a_matroid}}
We validate that $(\calV\setminus\calY, \calI')$ satisfies the conditions in Definition~\ref{def:matroid} of a matroid.  In particular:  
\begin{itemize}
\item to validate the first condition in Definition~\ref{def:matroid}, assume a set $\calX \subseteq \calV\setminus\calY$ such that $\calX\in \calI'$; moreover, assume a set $\calZ\subseteq \calX$; we need to show that $\calZ\in\calI'$\!\!. To this end, observe that the definition of $\calI'$ implies $\calX\cup\calY\in\calI$, since we assumed $\calX\in\calI'$\!\!. In~addition, the assumption $\calZ\subseteq \calX$ implies $\calZ\cup \calY\subseteq \calX\cup \calY$, and, as a result, $\calZ\cup \calY\in \calI$, since $(\calV,\calI)$ is a matroid.  Overall, $\calZ\subseteq \calV\setminus\calY$ (since $\calZ\subseteq \calX$, by assumption, and $\calX\subseteq \calV\setminus \calY$) and $\calZ\cup \calY\in \calI$; hence, $\calZ\in \calI'$\!\!, by the definition of $\calI'$\!\!, and now the first condition in Definition~\ref{def:matroid} is validated; 
\item to validate the second condition in Definition~\ref{def:matroid}, assume sets $\calX,\calZ\in\calV\setminus\calY$ such that $\calX,\calZ\in\calI'$ and $|\calX|<|\calZ|$; we need to show that there exists an element $z\in\calZ\setminus\calX$ such that $\calX\cup\{z\}\in\calI'$\!\!. To this end, observe that since $\calX,\calZ\in\calI'$\!\!, the definition of $\calI'$ implies that $\calX\cup\calY,\calZ\cup\calY\in\calI$.  Moreover, since $|\calX|<|\calZ|$, it also is  $|\calX\cup\calY|<|\calZ\cup\calY|$.  Therefore, since $(\calV,\calI)$ is a matroid, there exists an element $z\in (\calZ\cup\calY)\setminus(\calX\cup\calY)=\calZ\setminus\calX$ such that $(\calX\cup\calY)\cup\{z\}\in\calI$;
as a result, $\calX\cup\{z\}\in \calI'$\!\!, by the definition of $\calI'$\!\!.  
In sum,  $z\in \calZ\setminus\calX$ and $\calX\cup\{z\}\in \calI'$\!\!, and the second condition in Definition~\ref{def:matroid} is validated too.  

\hfill $\blacksquare$
\end{itemize}

\begin{mylemma}\label{lem:greedy_perfromance}
Recall the notation in Algorithm~\ref{alg:rob_sub_max}, and consider the sets $\calA_1$ and $\calA_2$ constructed by Algorithm~\ref{alg:rob_sub_max}'s lines~\ref{line:begin_while_1}-\ref{line:end_while_1} and lines~\ref{line:begin_while_2}-\ref{line:end_while_2}, respectively.  Then, for the set $\calA_2$, it holds:
\begin{itemize}
\item if the function $f$ is non-decreasing submodular and:
\begin{itemize}
\item if $(\calV,\calI)$ is a uniform matroid, then:
\begin{equation}\label{eq:greedy_perfromance_uniform}
f(\calA_{2})\geq \frac{1}{\kappa_f}(1-e^{-\kappa_f})\;\;\underset{\calX\subseteq \calV\setminus\calA_1, \calX\cup\calA_1\in \calI}{\max}  \; \;\; f(\calX).
\end{equation}
\item if $(\calV,\calI)$ is a matroid, then:
\begin{equation}\label{eq:greedy_perfromance_matroid}
f(\calA_{2})\geq \frac{1}{1+\kappa_f}\;\;\underset{\calX\subseteq \calV\setminus\calA_1, \calX\cup\calA_1\in \calI}{\max}  \; \;\; f(\calX).
\end{equation}
\end{itemize}
\item if the function $f$ is non-decreasing, then:
\begin{equation}\label{eq:greedy_perfromance}
f(\calA_{2})\geq (1-c_f)\;\;\underset{\calX\subseteq \calV\setminus\calA_1, \calX\cup\calA_1\in \calI}{\max}  \; \;\; f(\calX).
\end{equation}
\end{itemize}
\end{mylemma}
\paragraph*{Proof of Lemma~\ref{lem:greedy_perfromance}}
We first prove ineq.~\eqref{eq:greedy_perfromance}, then ineq.~\eqref{eq:greedy_perfromance_matroid}, and, finally, ineq.~\eqref{eq:greedy_perfromance_uniform}.
In particular, Algorithm~\ref{alg:rob_sub_max} constructs the set~$\calA_2$ greedily, by replicating the steps of the greedy algorithm introduced~\cite[Section~2]{fisher1978analysis}, to solve the following optimization problem: 
\begin{equation}
\label{eq:auxxxx}
\underset{\calX\subseteq \calV\setminus\calA_1, \calX\cup\calA_1\in \calI}{\max}  \; \;\; f(\calX);
\end{equation} 
let in the latter problem $\calI'\triangleq\{\calX\subseteq \calV\setminus\calA_1, \calX\cup\calA_1\in \calI\}$.  Lemma~\ref{lem:its_a_matroid} implies that $(\calA_1, \calI')$ is a matroid, and, as a result, the previous optimization problem is a matroid-constrained set function maximization problem. Now, to prove ineq.~\eqref{eq:greedy_perfromance}, ineq.~\eqref{eq:greedy_perfromance_matroid}, and ineq.~\eqref{eq:greedy_perfromance_uniform}, we make the following observations, respectively: when the function $f$ is merely non-decreasing, then \cite[Theorem~8.1]{sviridenko2017optimal} implies that the greedy algorithm introduced in~\cite[Section~2]{fisher1978analysis} returns for the optimization problem in eq.~\eqref{eq:auxxxx} a solution $\calS$ such that $f(\calS)\geq (1-c_f)\underset{\calX\subseteq \calV\setminus\calA_1, \calX\cup\calA_1\in \calI}{\max}\; f(\calX)$; this proves ineq.~\eqref{eq:greedy_perfromance}. Similarly, when the function $f$ is non-decreasing and submodular, then~\cite[Theorem~2.3]{conforti1984curvature} implies that the greedy algorithm introduced in~\cite[Section~2]{fisher1978analysis} returns for the optimization problem in eq.~\eqref{eq:auxxxx} a solution $\calS$ such that $f(\calS)\geq 1/(1+\kappa_f)\underset{\calX\subseteq \calV\setminus\calA_1, \calX\cup\calA_1\in \calI}{\max}\; f(\calX)$; this proves ineq.~\eqref{eq:greedy_perfromance_matroid}.  Finally, when the objective function $f$ is non-decreasing submodular, and when $\calI$ is a uniform matroid, then~\cite[Theorem~5.4]{conforti1984curvature} implies that the greedy algorithm introduced in~\cite[Section~2]{fisher1978analysis} returns for the optimization problem in eq.~\eqref{eq:auxxxx} a solution~$\calS$ such that $f(\calS)\geq 1/\kappa_f(1-e^{-\kappa_f})\underset{\calX\subseteq \calV\setminus\calA_1, \calX\cup\calA_1\in \calI}{\max} \; f(\calX)$; this proves ineq.~\eqref{eq:greedy_perfromance_uniform}, and concludes the proof of the lemma.
\hfill $\blacksquare$

\begin{mylemma}\label{lem:from_max_to_minmax}
Recall the notation in Theorem~\ref{th:alg_rob_sub_max_performance} and Appendix~\ref{app:notation}. Also, consider a uniform or partition matroid $(\calV,\calI')$. Then, for any set ${\calY}\subseteq \calV$ such that $\calY\in \calI$ and $\calY\in \calI'$\!\!, it holds:
\begin{equation}\label{eq:toprovefrom_max_to_minmax}
\underset{{\mathcal{X}}\subseteq \calV\setminus\calY, \calX\cup{{\mathcal{Y}}}\in \calI}{\max}  \; \;\; f({\mathcal{X}})\geq f(\calA^\star\setminus \calB^\star(\calA^\star)).
\end{equation}
\end{mylemma}
\paragraph*{Proof of Lemma~\ref{lem:from_max_to_minmax}}
We start from the left-hand-side of ineq.~\eqref{eq:toprovefrom_max_to_minmax}, and make the following observations:
\begin{align}
\underset{{\mathcal{X}}\subseteq \calV\setminus\calY, \calX\cup{{\mathcal{Y}}}\in \calI}{\max}  \;\; \; f({\mathcal{X}})&\geq \min_{\bar{\calY} \subseteq \calV, \bar{\calY}\in \calI,\calI'}\;\; \underset{\;{\mathcal{X}}\subseteq \calV\setminus\bar{\calY}, \calX\cup{\bar{\calY}}\in \calI}{\max}  \;\;\; f({\mathcal{X}})\nonumber\\
&=\min_{\bar{\calY} \subseteq \calV, \bar{\calY}\in \calI,\calI'}\;\; \underset{\bar{\mathcal{A}}\subseteq \calV, \bar{\mathcal{A}}\in \calI}{\max}  \;\;\; f({\bar{\mathcal{A}}\setminus \bar{\calY}})\nonumber\\
&\triangleq h.\nonumber
\end{align}
We next complete the proof of Lemma~\ref{lem:from_max_to_minmax} by proving that $h\geq f(\calA^\star\setminus \calB^\star(\calA^\star))$.  To this end, observe that for any set ${\calA}\subseteq \calV$ such that ${\calA}\in \calI$, and for any set ${\calY}\subseteq \calV$ such that $\calY\in \calI$ and $\calY\in \calI'$\!\!, it holds:
\begin{align}
\underset{\bar{\mathcal{A}}\subseteq \calV, \bar{\mathcal{A}}\in \calI}{\max}  \;\;\; f({\bar{\mathcal{A}}\setminus {\calY}})&\geq f({\mathcal{A}}\setminus {\calY}),\nonumber
\end{align}
which implies the following observations:
\begin{align}
h&\geq \min_{\bar{\calY} \subseteq \calV, \bar{\calY}\in \calI,\calI'}\;\; f({\mathcal{A}}\setminus \bar{\calY})\nonumber\\
&\geq \min_{\bar{\calY} \subseteq \calV, \bar{\calY}\in \calI'}\;\; f({\mathcal{A}}\setminus \bar{\calY})\nonumber\\
&=\min_{\bar{\calY} \subseteq \calA, \bar{\calY}\in \calI'}\;\; f({\mathcal{A}}\setminus \bar{\calY}),\nonumber
\end{align}
and, as a result, it holds:
\belowdisplayskip=-11pt\begin{align}
h&\geq \underset{\bar{\mathcal{A}}\subseteq \calV, \bar{\mathcal{A}}\in \calI}{\max}  \;\;\min_{\bar{\calY} \subseteq \calA, \bar{\calY}\in \calI'}\;\; f(\bar{\mathcal{A}}\setminus \bar{\calY})\nonumber\\
&= f(\calA^\star\setminus \calB^\star(\calA^\star)).\nonumber
\end{align}
\hfill $\blacksquare$

\section{Proof of Theorem~\ref{th:alg_rob_sub_max_performance}}

We first prove Theorem~\ref{th:alg_rob_sub_max_performance}'s part 1 (approximation performance), and then, Theorem~\ref{th:alg_rob_sub_max_performance}'s part 2 (running time).

\subsection{Proof of Theorem~\ref{th:alg_rob_sub_max_performance}'s part 1 (approximation performance)}

We first prove ineq.~\eqref{ineq:bound_non_sub}, and, then, ineq.~\eqref{ineq:bound_sub} and  ineq.~\eqref{ineq:bound_sub_uniform}.

To the above ends, we use the following notation (along with the notation in Algorithm~\ref{alg:rob_sub_max}, Theorem~\ref{th:alg_rob_sub_max_performance}, and Appendix~\ref{app:notation}):
\begin{itemize}
\item let $\calA_{1}^+\triangleq \calA_{1}\setminus \calB^\star(\calA)$, i.e., $\calA_{1}^+$ is the set of remaining elements in the set $\calA_{1}$ after the removal from $\calA_{1}$ of the elements in the optimal (worst-case) removal $\calB^\star(\calA)$;
\item let $\calA_{2}^+\triangleq \calA_{2}\setminus \calB^\star(\calA)$, i.e., $\calA_{2}^+$ is the set of remaining elements in the set $\calA_{2}$ after the removal from $\calA_{2}$ of the elements in the optimal (worst-case) removal $\calB^\star(\calA)$.
\end{itemize}

\bigskip

\paragraph*{Proof of ineq.~\eqref{ineq:bound_non_sub}}  Consider that the objective function~$f$ is non-decreasing and such that (without loss of generality) $f$ is non-negative and $f(\emptyset)=0$. Then, the proof of ineq.~\eqref{ineq:bound_non_sub} follows by making the following observations:
\belowdisplayskip=8pt\begin{align}
&\!\!\!f(\calA\setminus \calB^\star(\calA))\nonumber\\
&=f(\calA_{1}^+\cup \calA_{2}^+)\label{ineq2:aux_14}\\
&\geq (1-c_f)\sum_{v\in\calA_{1}^+\cup \calA_{2}^+}f(v)\label{ineq2:aux_15}\\
&\geq (1-c_f)\sum_{v\in\calA_{2}}f(v)\label{ineq2:aux_16}\\
&\geq (1-c_f)^2f(\calA_{2})\label{ineq2:aux_17}\\
&\geq (1-c_f)^3\underset{\calX\subseteq \calV\setminus\calA_1, \calX\cup\calA_1\in \calI}{\max}  \; \;\; f(\calX)\label{ineq2:aux_18}\\
&\geq (1-c_f)^3f(\calA^\star\setminus \calB^\star(\calA^\star)),\label{ineq2:aux_19}
\end{align}
where eqs.~\eqref{ineq2:aux_14}-\eqref{ineq2:aux_19} hold for the following reasons: eq.~\eqref{ineq2:aux_14} follows from the definitions of the sets~$\calA_{1}^+$ and $\calA_{2}^+$; ineq.~\eqref{ineq2:aux_15} follows from ineq.~\eqref{ineq2:aux_14} due to Lemma~\ref{lem:curvature}; ineq.~\eqref{ineq2:aux_16} follows from ineq.~\eqref{ineq2:aux_15} due to Lemma~\ref{lem:dominance}, which implies that for any element $v\in\calA_1^+$ and any element $v'\in \calA_2^+$ it is $f(v)\geq f(v')$
---note that due to the definitions of $\calA_{1}^+$ and $\calA_{2}^+$ it is $|\calA_{1}^+|=|\calA_{2}\setminus \calA_{2}^+|$, that is, the number of non-removed elements in~$\calA_{1}$ is equal to the number of removed elements in $\calA_{2}$,---
and the fact $\calA_{2}=(\calA_{2}\setminus \calA_{2}^+)\cup \calA_{2}^+$;
ineq.~\eqref{ineq2:aux_17} follows from ineq.~\eqref{ineq2:aux_16} due to Corollary~\ref{cor:ineq_from_lemmata}; ineq.~\eqref{ineq2:aux_18} follows from ineq.~\eqref{ineq2:aux_17} due to Lemma~\ref{lem:greedy_perfromance}'s  ineq.~\eqref{eq:greedy_perfromance}; finally, ineq.~\eqref{ineq2:aux_19} follows from ineq.~\eqref{ineq2:aux_18} due to Lemma~\ref{lem:from_max_to_minmax}. 
\hfill $\blacksquare$

\bigskip

In what follows, we first prove  ineq.~\eqref{ineq:bound_sub}, and then  ineq.~\eqref{ineq:bound_sub_uniform}: we first prove the part $\frac{1-\kappa_f}{1+\kappa_f}$ and $\frac{1-\kappa_f}{\kappa_f}(1-e^{-\kappa_f})$ of ineq.~\eqref{ineq:bound_sub} and of ineq.~\eqref{ineq:bound_sub_uniform}, respectively, and then, the part $\frac{h_f(\alpha,\beta)}{1+\kappa_f}$ and $\frac{h_f(\alpha,\beta)}{\kappa_f}(1-e^{-\kappa_f})$ of ineq.~\eqref{ineq:bound_sub} and  of ineq.~\eqref{ineq:bound_sub_uniform}, respectively.

\medskip

\paragraph*{Proof of part $(1-\kappa_f)/(1+\kappa_f)$ of ineq.~\eqref{ineq:bound_sub}}
Consider that the objective function~$f$ is non-decreasing submodular and such that (without loss of generality) $f$ is non-negative and $f(\emptyset)=0$.
To prove the part $(1-\kappa_f)/(1+\kappa_f)$ of
ineq.~\eqref{ineq:bound_sub} we follow similar observations to the ones we followed in the proof of ineq.~\eqref{ineq:bound_non_sub}; in particular:
\begin{align}
&\!\!\!f(\calA\setminus \calB^\star(\calA))\nonumber\\
&=f(\calA_{1}^+\cup \calA_{2}^+)\label{ineq5:aux_14}\\
&\geq (1-\kappa_f)\sum_{v\in\calA_{1}^+\cup \calA_{2}^+}f(v)\label{ineq5:aux_15}\\
&\geq (1-\kappa_f)\sum_{v\in\calA_{2}}f(v)\label{ineq5:aux_16}\\
&\geq (1-\kappa_f)f(\calA_{2})\label{ineq5:aux_17}\\
&\geq \frac{1-\kappa_f}{1+\kappa_f}\underset{\calX\subseteq \calV\setminus\calA_1, \calX\cup\calA_1\in \calI}{\max}  \; \;\; f(\calX)\label{ineq5:aux_18}\\
&\geq \frac{1-\kappa_f}{1+\kappa_f}f(\calA^\star\setminus \calB^\star(\calA^\star)),\label{ineq5:aux_19}
\end{align}
where eqs.~\eqref{ineq5:aux_14}-\eqref{ineq5:aux_19} hold for the following reasons: eq.~\eqref{ineq5:aux_14} follows from the definitions of the sets~$\calA_{1}^+$ and $\calA_{2}^+$; ineq.~\eqref{ineq5:aux_15} follows from ineq.~\eqref{ineq5:aux_14} due to Lemma~\ref{lem:non_total_curvature}; ineq.~\eqref{ineq5:aux_16} follows from ineq.~\eqref{ineq5:aux_15} due to Lemma~\ref{lem:dominance}, which implies that for any element $v\in\calA_1^+$ and any element $v'\in \calA_2^+$ it is $f(v)\geq f(v')$
---note that due to the definitions of the sets~$\calA_{1}^+$ and $\calA_{2}^+$ it is $|\calA_{1}^+|=|\calA_{2}\setminus \calA_{2}^+|$, that is, the number of non-removed elements in $\calA_{1}$ is equal to the number of removed elements in~$\calA_{2}$,--- and because $\calA_{2}=(\calA_{2}\setminus \calA_{2}^+)\cup \calA_{2}^+$; ineq.~\eqref{ineq5:aux_17} follows from ineq.~\eqref{ineq5:aux_16} because the set function $f$ is submodular, and as~a result, the~submodularity Definition~\ref{def:sub} implies that for any sets $\mathcal{S}\subseteq \mathcal{V}$ and $\mathcal{S}'\subseteq \mathcal{V}$, it is  $f(\mathcal{S})+f(\mathcal{S}')\geq f(\mathcal{S}\cup \mathcal{S}')$~\cite[Proposition 2.1]{nemhauser78analysis};  ineq.~\eqref{ineq5:aux_18} follows from ineq.~\eqref{ineq5:aux_17} due to Lemma~\ref{lem:greedy_perfromance}'s ineq.~\eqref{eq:greedy_perfromance_matroid}; finally, ineq.~\eqref{ineq5:aux_19} follows from ineq.~\eqref{ineq5:aux_18} due to Lemma~\ref{lem:from_max_to_minmax}.  
\hfill $\blacksquare$

\medskip

\paragraph*{Proof of part $(1-\kappa_f)/\kappa_f(1-e^{-\kappa_f})$ of  ineq.~\eqref{ineq:bound_sub_uniform}}
Consider that the objective function~$f$ is non-decreasing submodular and such that (without loss of generality) $f$ is non-negative and $f(\emptyset)=0$.  Moreover, consider that the pair $(\calV,\calI)$ is a uniform matroid. 
To prove the part $(1-\kappa_f)/\kappa_f(1-e^{-\kappa_f})$ of  ineq.~\eqref{ineq:bound_sub_uniform} we follow similar steps to the ones we followed in the proof of ineq.~\eqref{ineq:bound_sub} via the ineqs.~\eqref{ineq5:aux_14}-\eqref{ineq5:aux_19}.  We explain next where these steps differ: if instead of using Lemma~\ref{lem:greedy_perfromance}'s ineq.~\eqref{eq:greedy_perfromance_matroid} to get ineq.~\eqref{ineq5:aux_18} from ineq.~\eqref{ineq5:aux_17}, we use Lemma~\ref{lem:greedy_perfromance}'s ineq.~\eqref{eq:greedy_perfromance_uniform}, and afterwards apply Lemma~\ref{lem:from_max_to_minmax}, then, we derive ineq.~\eqref{ineq:bound_sub_uniform}.
\hfill $\blacksquare$

\medskip

\paragraph*{Proof of parts $h_f(\alpha,\beta)/(1+\kappa_f)$ and ${h_f(\alpha,\beta)}/{\kappa_f}(1-e^{-\kappa_f})$ of ineq.~\eqref{ineq:bound_sub} and  ineq.~\eqref{ineq:bound_sub_uniform}, respectively}

\begin{figure}[t]
\def \setAone{ (0,0) circle (1cm) }
\def \setBone{ (.5,0) circle (0.4cm)}
\def \setAtwo{ (2.5,0) circle (1cm) }
\def \setBtwo{ (3.0,0) circle (0.4cm)}
\def \myrectangle{ (-1.5, -1.5) rectangle (4, 1.5) }
\begin{center}
\begin{tikzpicture}
\draw \myrectangle node[below left]{$\mathcal{V}$};
\draw \setAone node[left]{$\mathcal{A}_{1}$};
\draw \setBone node[]{$\mathcal{B}_{1}^\star$};
\draw \setAtwo node[left]{$\mathcal{A}_{2}$};
\draw \setBtwo node[]{$\mathcal{B}_{2}^\star$};
\end{tikzpicture}
\end{center}
\caption{\small Venn diagram, where the sets $\mathcal{A}_{1}, \mathcal{A}_{2},\mathcal{B}_{1}^\star, \mathcal{B}_{2}^\star$ are as follows: per Algorithm~\ref{alg:rob_sub_max}, $\mathcal{A}_{1}$  and $\mathcal{A}_{2}$ are such that $\mathcal{A}=\mathcal{A}_{1}\cup \mathcal{A}_{2}$.  Due to their construction, it holds  $\mathcal{A}_{1}\cap \mathcal{A}_{2}=\emptyset$. Next, 
$\mathcal{B}_{1}^\star$ and $\mathcal{B}_{2}^\star$ are such that  $\mathcal{B}_{1}^\star=\mathcal{B}^\star(\mathcal{A})\cap\mathcal{A}_{1}$, and $\mathcal{B}_2^\star=\mathcal{B}^\star(\mathcal{A})\cap\mathcal{A}_{2}$; therefore, it is $\mathcal{B}_{1}^\star\cap \mathcal{B}_{2}^\star=\emptyset$ and $\mathcal{B}^\star(\mathcal{A})=(\mathcal{B}_{1}^\star\cup \mathcal{B}_{2}^\star)$.  
}\label{fig:venn_diagram_for_proof}
\end{figure}
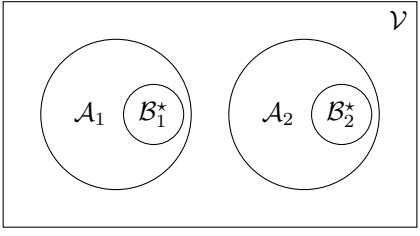

We complete the proof by first proving that:
\begin{equation}\label{ineq1:aux_1}
f(\calA\setminus\calB^\star(\calA))\geq \frac{1}{1+\beta}f(\calA_2),
\end{equation}
and, then, proving that:
\begin{equation}\label{ineq1:aux_2}
f(\calA\setminus\calB^\star(\calA))\geq \frac{1}{\alpha-\beta}f(\calA_2).
\end{equation}
The combination of ineq.~\eqref{ineq1:aux_1} and ineq.~\eqref{ineq1:aux_2} proves the part $\frac{h_f(\alpha,\beta)}{1+\kappa_f}$ and $\frac{h_f(\alpha,\beta)}{\kappa_f}(1-e^{-\kappa_f})$ of ineq.~\eqref{ineq:bound_sub} and  of ineq.~\eqref{ineq:bound_sub_uniform}, respectively, after also applying Lemma~\ref{lem:greedy_perfromance}'s ineq.~\eqref{eq:greedy_perfromance_matroid} and ineq.~\eqref{eq:greedy_perfromance_uniform}, respectively, and then Lemma~\ref{lem:from_max_to_minmax}.

\medskip

\textit{To prove ineq.~\eqref{ineq1:aux_1}}, we follow the steps of the proof of~\cite[Theorem~1]{tzoumas2017resilient}, and use the notation introduced in Fig.~\ref{fig:venn_diagram_for_proof}, along with the following notation:
\begin{align}
	\eta = \frac{f(\mathcal{B}_2^\star|\mathcal{A}\setminus \mathcal{B}^\star(\mathcal{A}))}{f(\mathcal{A}_2)}.
\end{align}

In particular, to prove ineq.~\eqref{ineq1:aux_1} we focus on the worst-case where $\calB^\star_2\neq \emptyset$; the reason is that if we assume otherwise, i.e.,  if we assume $\calB^\star_2= \emptyset$, then $f(\calA\setminus\calB^\star(\calA))=f(\calA_2)$, which is a tighter inequality to ineq.~\eqref{ineq1:aux_1}. Hence, considering $\calB^\star_2\neq \emptyset$,
we prove ineq.~\eqref{ineq1:aux_1} by first observing that:
\begin{equation}\label{ineq:aux_1}
f(\mathcal{A}\setminus\mathcal{B}^\star(\mathcal{A}))\geq\max\{f(\mathcal{A}\setminus\mathcal{B}^\star(\mathcal{A})),f(\mathcal{A}_1^+)\},
\end{equation}
and then proving the following three inequalities:
\begin{align}
f(\mathcal{A}\setminus\mathcal{B}^\star(\mathcal{A}))&\geq(1-\eta)f(\mathcal{A}_2)\label{ineq:aux_2},\\
f(\mathcal{A}_1^+)&\geq \eta \frac{1}{\beta}f(\mathcal{A}_2),\label{ineq:aux_3}\\
\max\{(1-\eta),\eta\frac{1}{\beta}\}&\geq \frac{1}{\beta+1}.\label{ineq:aux_4}
\end{align}

Specifically, if we substitute ineqs.~\eqref{ineq:aux_2}-\eqref{ineq:aux_4} to ineq.~~\eqref{ineq:aux_1}, and take into account that $f(\mathcal{A}_2)\geq 0$, then:
\begin{equation*}
f(\mathcal{A}\setminus\mathcal{B}^\star(\mathcal{A}))\geq \frac{1}{\beta+1}f(\mathcal{A}_2),
\end{equation*}
which implies ineq.~\eqref{ineq1:aux_1}.

We complete the proof of ineq.~\eqref{ineq:aux_1} by proving $0\leq \eta\leq 1$, and ineqs.~\eqref{ineq:aux_2}-\eqref{ineq:aux_4}, respectively.

\paragraph{Proof of ineq.~$0\leq \eta\leq 1$} We first prove that $\eta\geq 0$, and then, that $\eta\leq 1$: it holds~$\eta\geq 0$, since by definition $\eta=f(\mathcal{B}_2^\star|\mathcal{A}\setminus \mathcal{B}^\star(\mathcal{A}))/f(\mathcal{A}_2)$, and since $f$ is non-negative; and it holds~$\eta\leq 1$, since $f(\mathcal{A}_2)\geq f(\mathcal{B}^\star_2)$, due to monotonicity of $f$ and that $\mathcal{B}^\star_2 \subseteq \mathcal{A}_2$, and since $f(\mathcal{B}^\star_2)\geq f(\mathcal{B}_2^\star|\mathcal{A}\setminus \mathcal{B}^\star(\mathcal{A}))$, due to submodularity of $f$ and that $\emptyset \subseteq \mathcal{A}\setminus \mathcal{B}^\star(\mathcal{A})$. 

\paragraph{Proof of ineq.~\eqref{ineq:aux_2}}  We complete the proof of ineq.~\eqref{ineq:aux_2} in two steps.  First, it can be verified that:
\begin{align}\label{eq:aux_1}
& f(\mathcal{A}\setminus\mathcal{B}^\star(\mathcal{A}))=f(\mathcal{A}_2)-\nonumber\\ & f(\mathcal{B}^\star_2|\mathcal{A}\setminus\mathcal{B}^\star(\mathcal{A}))+f(\mathcal{A}_1|\mathcal{A}_2)-f(\mathcal{B}^\star_1|\mathcal{A}\setminus\mathcal{B}^\star_1),
\end{align}
since for any $\mathcal{X}\subseteq \mathcal{V}$ and $\mathcal{Y}\subseteq \mathcal{V}$, it holds $f(\mathcal{X}|\mathcal{Y})=f(\mathcal{X}\cup \mathcal{Y})-f(\mathcal{Y})$. Second, eq.~\eqref{eq:aux_1} implies ineq.~\eqref{ineq:aux_2}, since $f(\mathcal{B}^\star_2|\mathcal{A}\setminus\mathcal{B}^\star(\mathcal{A}))=\eta f(\mathcal{A}_2)$, and $f(\mathcal{A}_1|\mathcal{A}_2)-f(\mathcal{B}^\star_1|\mathcal{A}\setminus\mathcal{B}^\star_1)\geq 0$.
The latter is true due to the following two observations:~$f(\mathcal{A}_1|\mathcal{A}_2)\geq f(\mathcal{B}_1^\star|\mathcal{A}_2)$, since $f$ is monotone and $\mathcal{B}_1^\star \subseteq \mathcal{A}_1$; and~$f(\mathcal{B}_1^\star|\mathcal{A}_2)\geq f(\mathcal{B}^\star_1|\mathcal{A}\setminus\mathcal{B}^\star_1)$, since $f$ is submodular and $\mathcal{A}_2\subseteq \mathcal{A}\setminus\mathcal{B}^\star_1$ (see also Fig.~\ref{fig:venn_diagram_for_proof}).

\paragraph{Proof of ineq.~\eqref{ineq:aux_3}} 
Since it is  $\calB^\star_2\neq \emptyset$ (and as a result, it also is $\calA^+_1\neq \emptyset$), and since for all elements $a \in \calA^+_1$ and all elements $b\in \calB^\star_2$ it is $f(a)\geq f(b)$, from Lemma~\ref{lem:D3} we have:
\begin{align}
f(\calB^\star_2|\calA^+_1)&\leq |\calB^\star_2|f(\calA^+_1)\nonumber\\
&\leq \beta f(\calA^+_1),\label{aux:111}
\end{align}
since $|\calB^\star_2|\leq \beta$.  Overall,
\begin{align}
f(\calA^+_1)&\geq \frac{1}{\beta}f(\calB^\star_2|\calA^+_1)\label{aux5:1}\\
&\geq \frac{1}{\beta}f(\calB^\star_2|\calA^+_1\cup \calA^+_2)\label{aux5:2}\\
&=\frac{1}{\beta}f(\mathcal{B}_2^\star|\mathcal{A}\setminus \mathcal{B}^\star(\mathcal{A}))\label{aux5:3}\\
&=\eta\frac{1}{\beta}f(\calA_2),\label{aux5:4}
\end{align}
where ineqs.~\eqref{aux5:1}-\eqref{aux5:4} hold for the following reasons: ineq.~\eqref{aux5:1} follows from ineq.~\eqref{aux:111}; ineq.~\eqref{aux5:2} holds since $f$ is submodular and $\calA_1^+\subseteq \calA_1^+\cup \calA_2^+$; eq.~\eqref{aux5:3} holds due to the definitions of the sets $\calA_1^+$, $\calA_2^+$ and $\mathcal{B}^\star(\mathcal{A})$; finally, eq.~\eqref{aux5:4} holds due to the definition of $\eta$. 

\paragraph{Proof of ineq.~\eqref{ineq:aux_4}}  Let $b=1/\beta$.  We complete the proof first for the case where 
$(1-\eta)\geq \eta b$, and then for the case $(1-\eta)<\eta b$: when $(1-\eta)\geq \eta b$, $\max\{(1-\eta),\eta b\}= 1-\eta$ and $\eta \leq 1/(1+b)$; due to the latter, $1-\eta \geq b/(1+b)=1/(\beta+1)$ and, as a result,~\eqref{ineq:aux_4} holds. Finally, when $(1-\eta)< \eta b$, $\max\{(1-\eta),\eta b\}= \eta b$ and $\eta > 1/(1+b)$; due to the latter, $\eta b >  b/(1+b)$ and, as a result,~\eqref{ineq:aux_4} holds.

We completed the proof of~$0\leq \eta\leq 1$, and of ineqs.~\eqref{ineq:aux_2}-\eqref{ineq:aux_4}.  Thus, we also completed the proof of ineq.~\eqref{ineq1:aux_1}.

\medskip

\textit{To prove ineq.~\eqref{ineq1:aux_2}}, we consider the following mutually exclusive and collectively exhaustive cases:
\begin{itemize}
\item consider $\calB^\star_2= \emptyset$, i.e., all elements in $\calA_1$ are removed, and as result, none of the elements in $\calA_2$ is removed.  Then, $f(\calA\setminus \calB^\star(\calA))=f(\calA_2)$, and ineq.~\eqref{ineq1:aux_2} holds.
\item Consider $\calB^\star_2\neq \emptyset$, i.e., at least one of the elements in $\calA_1$ is \textit{not} removed; call any of these elements $s$.  Then:
\begin{equation}\label{eq1:aux1}
f(\calA\setminus\calB^\star(\calA))\geq f(s),
\end{equation}
since $f$ is non-decreasing. 
In addition:
\begin{equation}\label{eq1:aux2}
f(\calA_2)\leq \sum_{v\in \calA_2}f(v)\leq (\alpha-\beta) f(s),
\end{equation}
where the first inequality holds since $f$ is submodular~\cite[Proposition~2.1]{nemhauser78analysis}, and the second holds due to Lemma~\ref{lem:dominance} and the fact that $\calA_2$ is constructed by Algorithm~\ref{alg:rob_sub_max} such that $\calA_1\cup\calA_2\subseteq \calV$ and $\calA_1\cup\calA_2\in\calI$, where $|\calA_1|=\beta$ (since $\calA_1$ is constructed by Algorithm~\ref{alg:rob_sub_max} such that $\calA_1\subseteq \calV$ and $\calA_1\in\calI'$, where $(\calV,\calI')$ is a matroid with rank $\beta$) and $(\calV,\calI)$ is a matroid that has rank $\alpha$; the combination of ineq.~\eqref{eq1:aux1} and ineq.~\eqref{eq1:aux2} implies ineq.~\eqref{ineq1:aux_2}. 
\end{itemize}

Overall, the proof of ineq.~\eqref{ineq1:aux_2} is complete. \hfill $\blacksquare$

\subsection{Proof of Theorem~\ref{th:alg_rob_sub_max_performance}'s part 2 (running time)}

We complete the proof in two steps, where we denote the time for each evaluation of the objective function $f$ as $\tau_f$.  In particular, we first compute the running time of lines~\ref{line:begin_while_1}-\ref{line:end_while_1} and, then, of lines~\ref{line:begin_while_2}-\ref{line:end_while_2}: lines~\ref{line:begin_while_1}-\ref{line:end_while_1} need at most $|\calV|[|\calV|\tau_f+|\calV|\log(|\calV|)+|\calV|+O(\log(|\calV|))]$ time, since they are repeated at most $|\calV|$ times, and at each repetition line~\ref{line:select_element_bait} asks for at most $|\calV|$ evaluations of $f$\!, and for their sorting, which takes $|\calV|\log(|\calV|)+|\calV|+O(\log(|\calV|))$ time, using, e.g., the merge sort algorithm.  Similarly, lines~\ref{line:begin_while_2}-\ref{line:end_while_2} need $|\calV|[|\calV|\tau_f+|\calV|\log(|\calV|)+|\calV|+O(\log(|\calV|))]$.
Overall, Algorithm~\ref{alg:rob_sub_max} runs in $2|\calV|[|\calV|\tau_f+|\calV|\log(|\calV|)+|\calV|+O(\log(|\calV|))]=O(|\calV|^2\tau_f)$ time. 
\hfill $\blacksquare$


\begin{biography}[{\includegraphics[width=1in,height=1.25in,keepaspectratio]{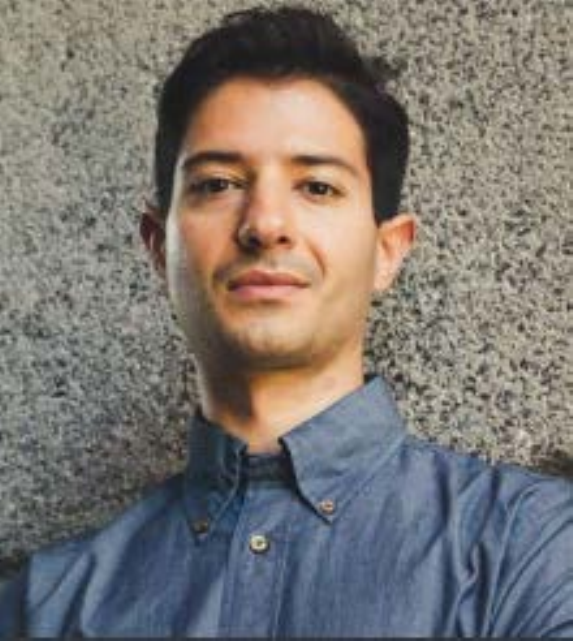}}]{Vasileios Tzoumas} (S'12-M'18) is a post-doctoral associate at the MIT Laboratory for Information \& Decision Systems (LIDS). He received his Ph.D.~at the Department of Electrical and Systems Engineering, University of Pennsylvania (2018). From July 2017 to December 2017 he was a visiting Ph.D.~student at the Institute for Data, Systems, and Society, MIT.  He holds a diploma in Electrical and Computer Engineering from the National Technical University of Athens (2012); a Master of science in Electrical Engineering from the University of Pennsylvania (2016); and a master of arts in Statistics from the Wharton School of Business at the University of Pennsylvania (2016).
Vasileios was a Best Student Paper Award finalist in the 56th IEEE Conference in Decision and Control (2017). 
\end{biography}

\begin{biography}[{\includegraphics[width=1in,height=1.25in,clip,keepaspectratio]{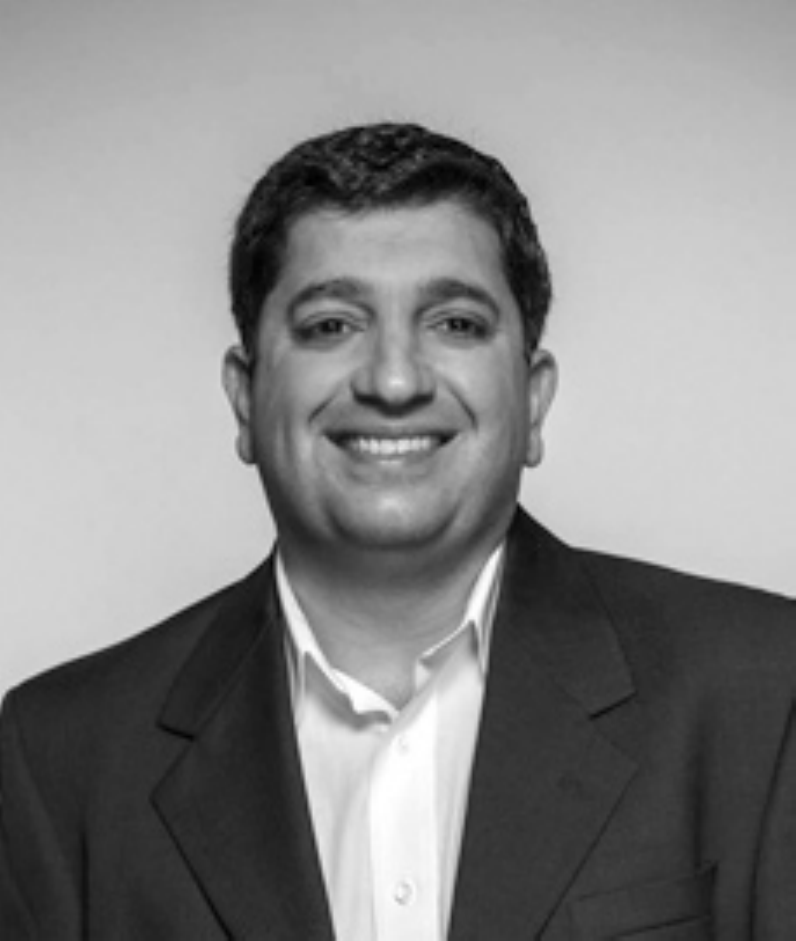}}] {Ali Jadbabaie} (S'99-M'08-SM'13-F'15) is the JR East Professor of Engineering and Associate Director of the Institute for Data, Systems and Society at MIT, where he is also on the faculty of the department of civil and environmental engineering and a principal investigator in the Laboratory for Information and Decision Systems (LIDS). He is the director of the Sociotechnical Systems Research Center, one of MIT's 13 laboratories. He received his Bachelors (with high honors) from Sharif University of Technology in Tehran, Iran, a Masters degree in electrical and computer engineering from the University of New Mexico, and his Ph.D.~in control and dynamical systems from the California Institute of Technology. He was a postdoctoral scholar at Yale University before joining the faculty at Penn in July 2002. Prior to joining MIT faculty, he was the Alfred Fitler Moore a Professor of Network Science and held secondary appointments in computer and information science and operations, information and decisions in the Wharton School. He was the inaugural editor-in-chief of IEEE Transactions on Network Science and Engineering, a new interdisciplinary journal sponsored by several IEEE societies. He is a recipient of a National Science Foundation Career Award, an Office of Naval Research Young Investigator Award, the O. Hugo Schuck Best Paper Award from the American Automatic Control Council, and the George S. Axelby Best Paper Award from the IEEE Control Systems Society. His students have been winners and finalists of student best paper awards at various ACC and CDC conferences. He is an IEEE fellow and a recipient of the Vannevar Bush Fellowship from the office of Secretary of Defense. His current research interests include the interplay of dynamic systems and networks with specific emphasis on multi-agent coordination and control, distributed optimization, network science, and network economics.
\end{biography}

\begin{biography}[{\includegraphics[width=1in,height=1.25in,clip,keepaspectratio]{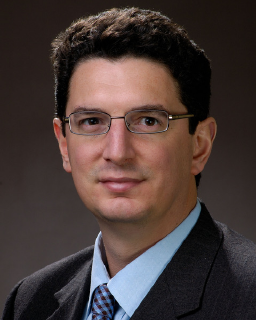}}]{George J.~Pappas} (S'90-M'91-SM'04-F'09) received the Ph.D.~degree in electrical engineering and computer sciences from the University of California, Berkeley, CA, USA, in  1998. He is currently the Joseph Moore Professor and Chair of the Department of Electrical and Systems Engineering, University of Pennsylvania, Philadelphia, PA, USA. He  also holds a secondary appointment with the Department of Computer and Information Sciences and the Department of Mechanical Engineering and Applied Mechanics. He is a Member of the GRASP Lab and the PRECISE Center. He had previously served as the Deputy Dean for Research with the School of Engineering and Applied Science. His research interests include control theory and, in particular, hybrid systems, embedded systems, cyberphysical systems, and hierarchical and distributed control systems, with applications to unmanned aerial vehicles, distributed robotics, green buildings, and biomolecular networks. Dr. Pappas has received various awards, such as the Antonio Ruberti Young Researcher Prize, the George S. Axelby Award, the Hugo Schuck Best Paper Award, the George H. Heilmeier Award, the National Science Foundation PECASE award and numerous best student papers awards.
\end{biography}

\end{document}